\documentclass[12pt,reqno]{amsart}
\usepackage{amscd,amssymb}

\swapnumbers
\newtheorem{thm}[subsection]{Theorem}
\newtheorem{cor}[subsection]{Corollary}
\newtheorem{lem}[subsection]{Lemma}
\newtheorem{prop}[subsection]{Proposition}

\theoremstyle{definition}

\numberwithin{equation}{subsection}
\newcommand{\tr}{\sideset{^{t{}}}{}\tnull\nolimits\hskip-3.5pt}

\DeclareMathOperator*{\tnull}{}
\renewcommand{\mod}[1]{ (\text{\rm mod }#1)}
\newcommand{\quash}[1]{}
\newcommand{\ua}{\frac{\textstyle{u_1}}{\textstyle{N}}}
\newcommand{\ub}{\frac{\textstyle{u_2}}{\textstyle{N}}}
\newcommand{\va}{\frac{\textstyle{v_1}}{\textstyle{N}}}
\newcommand{\vb}{\frac{\textstyle{v_2}}{\textstyle{N}}}
\renewcommand{\k}[1]{\kappa\left( #1 \right)}

\begin{document}

\title{Abelian surfaces with anti-holomorphic multiplication}
\author{Mark Goresky${}^1$}\thanks{1.  School of Mathematics, Institute for Advanced Study,
Princeton N.J.  Research partially supported by NSF grants \# DMS 9626616 and DMS
9900324.}
\author{Yung sheng Tai${}^2$}\thanks{2.  Dept. of Mathematics, Haverford College, Haverford
PA.}
\quash{
\begin{abstract}
For appropriate $N\ge 3$ and $d<0,$ the moduli space of principally
polarized abelian surfaces with level $N$ structure and anti-holomorphic multiplication by
$\mathcal O_d$ (the ring of integers in $\mathbb Q(\sqrt{d})$) is shown to consist of the real
points of a quasi-projective algebraic variety defined over $\mathbb Q$, and to coincide with
finitely many copies of the quotient of hyperbolic 3-space by the principal congruence
subgroup of level $N$ in $\mathbf{SL}(2, \mathcal O_d).$
\end{abstract}
} %% end quash
\maketitle

\section{Introduction}
Let $\mathfrak h_2 = \mathbf{Sp}(4,\mathbb R)/\mathbf{U}(2)$ be the Siegel upperhalf space of
rank $2$.  The quotient $\mathbf{Sp}(4,\mathbb Z) \backslash\mathfrak h_2$ has three
remarkable properties:  (a)  it is the moduli space of principally polarized abelian surfaces,
(b) it has the structure of a quasi-projective complex algebraic variety which is defined over
the rational numbers $\mathbb Q,$ and (c) it has a natural compactification (the Baily-Borel
Satake compactification) which is defined over the rational numbers. 

Fix a square-free integer $d < 0.$  One might ask whether similar statements hold for the
arithmetic quotient 
\begin{equation}\label{eqn-W}
W = \mathbf{SL}(2,\mathcal O_d) \backslash \mathbf{H_3} \end{equation}
where $\mathbf{H_3} = \mathbf{SL}(2,\mathbb C)/ \mathbf{SU}(2)$ is the hyperbolic 3-space and
where $\mathcal O_d$ is the ring of integers in the quadratic imaginary number field $\mathbb
Q(\sqrt{d}).$  One might first attempt to  interpret (\ref{eqn-W}) as a moduli space for
abelian surfaces with complex multiplication, but this is wrong.  Moreover, the space $W$ is
(real) 3-dimensional and it does not have an algebraic structure.  Nevertheless, in this paper
we show that there are analogues to all three statements (a), (b), and (c) if we introduce the
appropriate ``real'' structure and level structure on the abelian varieties.

Consider a 2-dimensional abelian variety $A$ with a principal polarization $H$ and with a
homomorphism $\psi : \mathcal O_d \to \text{End}_{\mathbb R}(A)$ such that $\kappa =
\psi(\sqrt{d})$ acts as an anti-holomorphic endomorphism of $A.$  Such a triple $(A, H,
\kappa)$ will be referred to (in \S \ref{sec-anticomplex}) as a principally polarized abelian
surface with anti-holomorphic multiplication by $\mathcal O_d$    There is an associated
notion of a level $N$ structure on such a triple.  (See \S \ref{sec-modulispace}.)

If $N$ is a positive integer, let
\[ \Lambda_N = \mathbf{SL}(2,\mathcal O_d)[N] \text{ and } \Gamma(N) = \mathbf{Sp}(4,\mathbb
Z)[N]\]
be the principal congruence subgroups of $\mathbf{SL}(2,\mathcal O_d)$ and
$\mathbf{Sp}(4,\mathbb Z)$ respectively, with level $N.$  Fix $d < 0$ and assume that $\mathbb
Q(\sqrt{d})$ has class number one.  Suppose $N \ge 3$ and, if $d \equiv 1 \mod 4,$ then assume
also that $N$ is even.  For simplicity (and for the purposes of this introduction only),
assume that $d$ is invertible $\mod N.$  In Theorems \ref{thm-realpoints} and \ref{thm-main}
we prove the following analogues of statements (a) and (b) above.

\begin{thm}  Given $d,N$ as above, the moduli space $V(d,N)$ of principally polarized abelian
surfaces with anti-holomorphic multiplication by $\mathcal O_d$ and level $N$ structure
consists of finitely many copies of the arithmetic quotient $\Lambda_N \backslash
\mathbf{H_3}.$  These copies are indexed by a certain (nonabelian cohomology) set $H^1(\mathbb
C/ \mathbb R, \Gamma(N)).$   Moreover this moduli space $V(d,N)$ coincides with the set of
real points $X_{\mathbb R}$ of a quasi-projective complex algebraic variety which is defined
over $\mathbb Q.$ 
\end{thm}

The key observation is that there is an involution $\ \widehat{}\ $ on $\mathbf{Sp}(4,\mathbb
R)$ whose fixed point set is $\mathbf{SL}(2,\mathbb C)$ and which has the property that the
fixed subgroup of $\Gamma(N)$ is exactly $\Lambda_N.$  This observation was inspired by
Nygaard's earlier article \cite{Nygaard}.  The key technical tool, Proposition
\ref{prop-Comessatti}, an analog of the lemma of Comessatti and Silhol (\cite{Silhol2}),
describes ``normal forms'' for the period matrix of an abelian surface with anti-holomorphic
multiplication.

The algebraic variety $X_{\mathbb C}$  is just the Siegel moduli space $\Gamma(N) \backslash
\mathfrak h_2$ of principally polarized abelian surfaces with level $N$ structure.  The
involution $\ \widehat{}\ $ extends to an anti-holomorphic involution of the Baily-Borel
Satake compactification $\overline{X}$ of $ X_{\mathbb C}$ and hence defines a real
structure on $\overline{X}.$  In \S \ref{sec-rational} we prove an analogue of statement (c)
above by showing that $\overline{X}$ admits a rational structure which is compatible with this
real structure.

One might ask whether the topological closure $\overline{V}$ of $V(d,N)$ in $\overline{X}$
coincides with the set of real points $\overline{X}_{\mathbb R}$ of the Baily-Borel
compactification.  At the moment, we do not know the answer to this question.  However, in
Theorem \ref{thm-noboundary} we show that the difference $\overline{X}_{\mathbb R} -
\overline{V}$ consists at most of finitely many isolated points.

The results in this article are, to a large extent, parallel to those of \cite{GT}, in which
similar phenomena are explored for $\mathbf{GL}(n, \mathbb R).$  However the techniques used
here are completely different from those in \cite{GT} and we have been unable to find a common
framework for both papers.  Nevertheless we believe they are two examples of a more general
phenomenon in which an arithmetic quotient of a (totally real) Jordan algebra admits the
structure of a (connected component of a) real algebraic variety which classifies real abelian
varieties with polarization, endomorphism, and level structures.

The authors would like to thank the Institute for Advanced Study for its support and
hospitality.

\section{An involution on the symplectic group}
\subsection{}\label{subsec-starthere}
We identify the symplectic group $\mathbf{Sp}(4,\mathbb R)$ with all 4 $\times$ 4 real
matrices $g=\left( \begin{smallmatrix} A & B \\ C & D \end{smallmatrix} \right)$ such that
\begin{equation}\label{eqn-sp4}
\tr{A}D - \tr{C}B = I;\ \tr{A}C \text{ and } \tr{B}D \text{ are symmetric}.
\end{equation}
The inverse of the symplectic matrix $g$ is $\left( \begin{smallmatrix}
\tr{D} & -\tr{B} \\ - \tr{C} & \tr{A} \end{smallmatrix} \right).$  The symplectic group acts
transitively on the Siegel upper halfspace
\begin{equation*}\mathfrak h_2 = \left\{Z=
X+iY \in M_{2\times 2}(\mathbb C)\ \Bigm|\ \tr Z=Z \text{ and } Y>0  \right\}
\end{equation*} 
by fractional linear transformations: 
\begin{equation*}
g\cdot Z = (AZ+B)(CZ+D)^{-1}.\end{equation*} 
The stabilizer of the basepoint 
\[ e_0 = iI \]
is the unitary group $\text{U}(2)$ which is embedded in $\mathbf{Sp}(4,\mathbb R)$ by
$A +iB \mapsto \left( \begin{smallmatrix} A & B \\ -B & A \end{smallmatrix} \right).$

\subsection{The number field}\label{subsec-numberfield}
Fix a square-free integer $d<0$ and let $\mathcal O_d$ be the ring of integers in the
quadratic imaginary number field $\mathbb Q(\sqrt d)$, that is,
\begin{equation*}
\mathcal O_d = \begin{cases}
 \mathbb Z[\sqrt d] &\text{ if } d \not\equiv 1 \mod 4 \\
\mathbb Z[\frac{1 + \sqrt d}{2}] &\text{ if } d \equiv 1 \mod 4 
\end{cases} \end{equation*}
Throughout the rest of this paper we will need to consider the two cases $d \not\equiv 1 \mod
4$ and $d \equiv 1 \mod 4$ separately.
\subsection{The embeddings}\label{subsec-embeddings}
If $d \not\equiv 1\mod 4$ define
$ \sigma= I_2$ to be the $2\times 2$ identity matrix and $M = I_4$ to be the $4\times 4$
identity matrix.  If $d \equiv 1 \mod 4$ define
$\sigma= \left( \begin{smallmatrix} 2 & 0 \\ -1 & 1 \end{smallmatrix} \right)$ and
\[ M = \left( \begin{matrix} \sigma & 0 \\ 0 & 2 \tr{\sigma}^{-1} \end{matrix} \right)\in
\mathbf{GL}(4,\mathbb R).\]  Let 
\[ \mathbf H_3 = \{ (z,r) \in \mathbb C \times \mathbb R| \ \text{Im}(z) > 0 \} \]
denote the hyperbolic 3-space.  Define the embedding $\phi: \mathbf H_3 \hookrightarrow
\mathfrak h_2$ by
\[ \phi(z,r) = \sigma \left( \begin{matrix} z & r \\ r & d\bar z \end{matrix} \right)
\tr{\sigma}. \]
Define the embedding $\phi:\mathbf{SL}(2,\mathbb C) \to \mathbf{Sp}(4,\mathbb R)$ by
\begin{equation}\label{eqn-phi} \phi \left( \begin{matrix}
a_1 + b_1\sqrt d & a_2+b_2\sqrt d
\\ a_3 + b_3\sqrt d & a_4 + b_4\sqrt d 
\end{matrix} \right) =  M\left(\begin{matrix}
a_1 & b_1 & b_2 & a_2\\ b_1d & a_1 & a_2 & b_2d\\ b_3d & a_3 & a_4 & b_4d\\
a_3 & b_3 & b_4 & a_4
\end{matrix} \right)M^{-1}. \end{equation}
The matrix $M$ normalizes $\mathbf{Sp}(4,\mathbb R)$ and was chosen so that Proposition
\ref{prop-arithmeticgroup} (below) holds.

\subsection{The involutions}\label{subsec-involutions}
If $d \not\equiv 1 \mod{4}$ define $\beta = \left( \begin{smallmatrix} 0 & 1 \\ d & 0
\end{smallmatrix} \right).$  If $d=4m+1$ define $\beta = \left( \begin{smallmatrix} 1 & 2 \\
2m & -1 \end{smallmatrix} \right).$  Define
\begin{equation}\label{eqn-Nbeta}
 N_{\beta} = \left( \begin{matrix} \beta & 0 \\ 
0 & \tr{\beta} \end{matrix} \right).\end{equation}
Then $\beta^2 = dI_2$ and $N_{\beta}^2 = dI_4$ and $\sigma^{-1} \beta \sigma = \left(
\begin{smallmatrix} 0 & 1 \\ d & 0 \end{smallmatrix} \right).$  The element $N_{\beta}$
normalizes $\mathbf{Sp}(4,\mathbb R)$ and gives rise to an involution
on $\mathbf{Sp}(4,\mathbb R)$,
\begin{equation}\label{eqn-theinvolution}
 \widehat{g} = N_{\beta} g N_{\beta}^{-1} = N_{\beta}^{-1} g N_{\beta}.  \end{equation}
It also gives rise to an involution on $\mathfrak h_2$,
\[ \widehat{Z} = \beta \overline{Z} \tr{\beta}^{-1}= \frac{1}{d}\beta\overline{Z}\tr{\beta}.\]

\begin{prop}\label{prop-embeddings} The embeddings (\ref{subsec-embeddings}) and involutions
(\ref{subsec-involutions}) are related as follows.
\begin{enumerate}
\item  For any $g \in \mathbf{Sp}(4,\mathbb R)$ and $x \in \mathfrak h_2$ we have:
$ \widehat{gx} = \widehat{g} \widehat{x}.$
\item  The fixed point set of the involution $\ \widehat{} \ $ on $\mathbf{Sp}(4,\mathbb R)$
is $\phi(\mathbf{SL}(2,\mathbb C)).$
\item  The fixed point set of the involution $\ \widehat{}\ $ on $\mathfrak h_2$ is
$\phi(\mathbf{H}_3).$
\item There is a unique (transitive) action of $\mathbf{SL}(2,\mathbb C)$ on $\mathbf H_3$ so
that $\phi(gx) = \phi(g) \cdot \phi(x)$ for all $g \in \mathbf{SL}(2,\mathbb C)$ and all $x
\in\mathbf H_3.$ \end{enumerate}\end{prop}

\subsection{Proof}
The proof is a direct calculation for which it helps to observe that
\[
M^{-1} N_{\beta} M = \left( \begin{array}{cc|cc}
0 & 1 \\ d & 0 \\ \hline & & 0 & d \\ & & 1 & 0 \end{array} \right). \qed
\]

\subsection{The basepoint}
If $h = \left( \begin{smallmatrix} X & 0 \\ 0 & Y \end{smallmatrix} \right) \in
\mathbf{GL}(4,\mathbb R)$ normalizes $\mathbf{Sp}(4,\mathbb R)$, if $ e \in \mathfrak h_2$ and
if $XeY^{-1} \in \mathfrak h_2$ then we will write $h\cdot e = XeY^{-1}$ and observe that
\[ \text{Stab} (XeY^{-1}) = h \text{Stab}(e) h^{-1}\]
where $\text{Stab}(e)$ denotes the stabilizer of $e$ in $\mathbf{Sp}(4,\mathbb R).$
Let $\tau = \left( \begin{smallmatrix} 1 & 0 \\ 0 & \sqrt{-d} \end{smallmatrix} \right)$ and
$T = \left( \begin{smallmatrix} \tau & 0 \\ 0 & \tau^{-1} \end{smallmatrix} \right) \in
\mathbf{Sp}(4,\mathbb R).$  Then $MT$ normalizes $\mathbf{Sp}(4,\mathbb R).$  Define the
basepoint
\[ e_1 = MT \cdot e_0\]
where $e_0 = iI$ as in \S \ref{subsec-starthere}.  Define $u_0 = \left( \begin{smallmatrix} 0
& 1 \\ -1 & 0 \end{smallmatrix} \right);$ with $\tr{u}_0 = u_0^{-1} = -u_0.$   The proof of
the following lemma (which will be needed in \S \ref{subsec-stab}) is a direct calculation.

\begin{lem}\label{lem-basepoint}
The basepoint $e_1$ satisfies $\widehat{e_1} = e_1.$  If $s = (MT)c(MT)^{-1} \in
\text{Stab}(e_1)$ with $c \in \text{U}(2)$ then
\[ \widehat{s} = (MT) \overset{\bullet}c (MT)^{-1} \]
where   $c \mapsto \overset{\bullet}c$ is the involution on $\text{U}(2)$ which is given by
\begin{equation}
\label{eqn-bullet} \overset{\bullet}c = u_0 \overline{c} u_0^{-1} = - u_0 \tr{c}^{-1} u_0.
\qed \end{equation}
\end{lem}

\subsection{Arithmetic subgroups}\label{subsec-arithmeticsubgroups}
Fix a (rational) integer $N \in \mathbb Z.$  Let $(N) \subset \mathcal O_d$ be the ideal
generated by $N.$  Let $z = a + b \sqrt{d} \in \mathcal O_d.$  If $d \not\equiv 1\mod{4}$ then
$a,b$ are integers, and \begin{align}
 z \equiv 1\mod{(N)} &\text{ iff } N|(a-1) \text{ and } N|b.\\
\intertext{If $d \equiv 1 \mod{4}$ then  $a,b$ are both integers or both half-integers, and}
\label{eqn-modN}
z \equiv 1\mod{(N)} &\text{ iff } N|(a-b-1) \text{ and } N|2b \end{align}
(which may be seen by writing $z = u + v \left( \frac{1 + \sqrt{d}}{2} \right)$ with $N|(u-1)$
and $N|v$).

Define the principal congruence subgroup
\begin{equation}\label{eqn-LambdaN} 
\Lambda_N=\mathbf{SL}(2,\mathcal O_d)(N) = \left\{\left( \begin{matrix} z_1 & z_2 \\
z_3 & z_4 \end{matrix}\right)\in \mathbf{SL}(2,\mathbb C) \bigg| z_i \in \mathcal O_d,
\begin{matrix} 
z_1, z_4 \equiv 1 \mod{(N)} \\
z_2, z_3 \equiv 0 \mod{(N)} \end{matrix} \right\}.
 \end{equation}
For $N=1$ set $\Lambda=\mathbf{SL}(2,\mathcal O_d)(N) = \mathbf{SL}(2,\mathcal O_d).$

Let $\Gamma(N) \subset \mathbf{Sp}(4,\mathbb Z)$ denote the principal congruence subgroup of
level $N.$ (For $N=1$ set $\Gamma(N) = \mathbf{Sp}(4,\mathbb Z).$)  Let $\widehat{\Gamma}(N) =
\{ \widehat{\gamma}|\ \gamma \in \Gamma(N)\}.$  Define 
\begin{equation}\label{eqn-GammaN}
 \Gamma_N = \Gamma(N) \cap \widehat{\Gamma}(N).\end{equation}
If $d$ is invertibile  $\mod N$ then $\Gamma_N = \Gamma(N).$  For, if $g \in \Gamma(N)$ then
\[d \widehat{g} =  N_{\beta} g N_{\beta} \equiv N_{\beta}^2 \equiv dI \mod N \]
hence $\widehat{g} \in \Gamma(N).$  In general, $\Gamma_N$ has finite index in $\Gamma(N)$
because $\Gamma(dN) \subset \Gamma_N$ has finite index in $\Gamma(N).$

\begin{prop}\label{prop-arithmeticgroup}  The embedding $\phi:\mathbf{SL}(2,\mathbb C) \to
\mathbf{Sp}(4,\mathbb R)$ of \S
\ref{subsec-embeddings} satisfies:
\[ \phi^{-1}(\Gamma_N) = \Lambda_N=\mathbf{SL}(2,\mathcal O_d)(N) \]
for any $N \ge 1.$
\end{prop}

\subsection{Proof}
Let
\[ h = \left( \begin{matrix} a_1+b_1\sqrt{d} & a_2+b_2\sqrt{d} \\
a_3+b_3\sqrt{d} & a_4+b_4 \sqrt{d} \end{matrix} \right)\in \mathbf{SL}(2,\mathbb C). \]
  Set 
\begin{equation}\label{eqn-4matrices}
A = \left( \begin{smallmatrix}  a_1 & b_1 \\ db_1& a_1 \end{smallmatrix}\right),\ 
B = \left( \begin{smallmatrix} b_2 & a_2 \\ a_2 & db_2 \end{smallmatrix} \right),\
C = \left( \begin{smallmatrix} db_3 & a_3 \\ a_3 & b_3 \end{smallmatrix} \right),\
D = \left( \begin{smallmatrix} a_4 & db_4 \\ b_4 & a_4 \end{smallmatrix} \right)
\end{equation}
First suppose that $d \not\equiv 1 \mod 4.$  Then $ \phi(h) = \left( \begin{smallmatrix} A & B
\\ C & D \end{smallmatrix} \right).$
It follows that $\phi(h) \in \Gamma(N)$ iff $h \in \mathbf{SL}(2,\mathcal O_d)(N).$  This
shows that $\phi(\mathbf{SL}(2,\mathbb C)) \cap \Gamma(N) = \phi(\mathbf{SL}(2,\mathcal
O_d)(N)).$   But the image of $\phi$ is fixed under the involution $ g \mapsto \widehat{g}$ so
this intersection is also contained in $\widehat{\Gamma}(N)$ as claimed.

Now suppose that $d \equiv 1 \mod 4.$  Then
\begin{equation}\label{eqn-phih} \phi(h) = \left( \begin{matrix} 
\sigma A \sigma^{-1} &  \frac{1}{2}\sigma B \tr{\sigma} \\
2\tr{\sigma}^{-1}C \sigma^{-1} & \tr{\sigma}^{-1} D \tr{\sigma} \end{matrix}
\right).\end{equation}
Moreover, $h \in \mathbf{SL}(2,\mathcal O_d)$ iff for each $i$
($1 \le i \le 4$), either $a_i,b_i$ are both integral or they are both half-integral.  
A simple calculation using (\ref{eqn-modN}) shows that $\phi(h) \in \Gamma(N) \iff h \in
\mathbf{SL}(2,\mathcal O_d)(N).$  As in the preceding paragraph, the image of $\phi$ is fixed
under the involution $g \mapsto \widehat{g},$ which implies that $\phi(\mathbf{SL}(2,\mathbb
C)) \cap \Gamma(N)$ is also contained in $\widehat{\Gamma}(N).$  \qed

\subsection{$\Gamma$-real points}
Let $\gamma \in Sp_4(\mathbb R)$ and define the locus of $(\gamma,\ \widehat{} \ )$-real
points,
\[ E_{\gamma} = \left\{ Z \in \mathfrak h_2|\ \widehat{Z}=\gamma \cdot Z \right\}. \]
Then $E_I = \phi(\mathbf H_3)$ is the fixed point set of $Z \mapsto \widehat Z.$
If $\Gamma \subset Sp_4(\mathbb R)$ is an arithmetic subgroup, define the
$(\Gamma,\ \widehat{}\ )$-real points of $\mathfrak h_2$ to be the union
$\cup_{\gamma}E_{\gamma}$ over all $\gamma \in \Gamma$. 

\section{Galois cohomology}\subsection{}
The involution $g \mapsto \widehat{g}$ may be considered to be an action of
$\text{Gal}(\mathbb C / \mathbb R) = \{ 1,\ \widehat{}\ \}$ on $\mathbf{Sp}(4,\mathbb R)$.
For any $\gamma \in \mathbf{Sp}(4,\mathbb R)$ let $f_{\gamma}:\text{Gal}(\mathbb C/ \mathbb R)
\to \mathbf{Sp}(4,\mathbb R)$ be the mapping $f_{\gamma}(1)=I$ and $f_{\gamma}(\ \widehat{}\ )
=\gamma.$  Then $f_{\gamma}$ is a 1-cocycle iff $\gamma \widehat{\gamma} =I$ and it is a
coboundary iff there exists $h \in \mathbf{Sp}(4,\mathbb R)$ such that $\gamma = \widehat{h}
h^{-1}.$ 

\begin{lem}\label{subsec-stab}  With respect to the above action, the (nonabelian) Galois
cohomology set $H^1(\mathbb C/\mathbb R, \text{Stab}(e_1))$ is trivial. \end{lem}

\subsection{Proof}
By Lemma \ref{lem-basepoint} it suffices to show that $H^1(\mathbb C/\mathbb R, \mathbf{U}
(2))=0$ where $\text{Gal}(\mathbb C/\mathbb R)$ acts on $\mathbf{U}(2)$ by $c \mapsto
\overset{\bullet}{c}.$  Assume $c \overset{\bullet}{c} = I.$  We must show there exists $w \in
\mathbf{U}(2)$ so that $c = \overset{\bullet}{w}w^{-1}.$

Embed $U(1) \hookrightarrow \mathbf{U}(2)$ by  $e^{i\theta} \mapsto \left(\begin{smallmatrix}
e^{i\theta} & 0 \\ 0 & e^{i\theta} \end{smallmatrix}\right).$   From (\ref{eqn-bullet}), 
$c \overset{\bullet}{c} = I$ if and only if $ cu_0 = u_0 \tr{c}$ which holds iff $c \in U(1).$ 
(Write $c = \left(\begin{smallmatrix}
c_1& c_2\\ c_3 & c_4 \end{smallmatrix} \right)$ and compute both sides to get
$c_1 = c_4$ and $c_2 = c_3 =0.$)  If $c = e^{i\theta}$ then it suffices to take
$w = e^{-i\theta/2}.$  \qed

\begin{lem}\label{lem-coboundary}  Let $\gamma \in \mathbf{Sp}(4,\mathbb R)$ and suppose that
$f_{\gamma}$ is a 1-cocycle.  Then $f_{\gamma}$ is a coboundary iff $E_{\gamma} \ne \phi.$
\end{lem}

\subsection{Proof}   If $f_{\gamma}$ is a coboundary, say $\gamma = \widehat{w} w^{-1}$
then $E_{\gamma} \supset w\phi(\mathbf H_3)$ is not empty.  On the other hand, suppose
$\widehat{Z} = \gamma Z.$  Choose $k \in \mathbf{Sp}(4,\mathbb R)$ so that $Z = ke_1$ (hence
$\widehat{Z} = \widehat{k}e_1$).  Therefore $u:=\widehat{k}^{-1} \gamma k \in \text{Stab}
(e_1)$ and $u\widehat{u} = I.$  By \S \ref{subsec-stab} there exists $w \in
\mathbf{Sp}(4,\mathbb R)$ so that $u = \widehat{w} w^{-1}$, hence $\gamma =
\widehat{kw}(kw)^{-1}.$  \qed

\begin{prop}\label{prop-Galois}
  With respect to the above action, the (nonabelian) Galois cohomology set $H^1(\mathbb C/
\mathbb R, \mathbf{Sp}_4(\mathbb R))$ is trivial. \end{prop}

\subsection{Proof} Suppose $\gamma = \left(\begin{smallmatrix} A&B \\ C&D \end{smallmatrix}
\right) \in \mathbf{Sp}(4,\mathbb R)$ and that $\gamma \widehat{\gamma} = I.$  Then
\begin{equation}\label{eqn-conditions}
\beta A \beta^{-1} = \tr{D};\ \beta D \text{ and } \tr{\beta}C \text{ are skew symmetric}.
\end{equation}
Write $\tr{\beta}C = \left( \begin{smallmatrix} 0 & \lambda \\ -\lambda & 0 \end{smallmatrix}
\right).$  We consider two cases:  $\lambda \ne 0$ and $\lambda = 0.$

First suppose that $\lambda \ne 0.$  Then $C$ is invertible and $\left( \tr{\beta}
C\right)^{-1}  = - \frac{1}{\lambda^2} \tr{\beta}C.$  Set
\begin{equation*}\gamma' = \widehat h \gamma h^{-1}\text{ where }
h = \left( \begin{matrix} I & -C^{-1} D \\ 0 & I \end{matrix} \right). 
\end{equation*}
Then $f_{\gamma'}$ is also a 1-cocycle which is cohomologous to $f_{\gamma}.$  In fact,
\begin{equation*}
\gamma ' = \left( \begin{matrix} 0 & - \tr{C}^{-1} \\ C & 0 \end{matrix} \right)
\end{equation*}
which is most easily seen by writing $\gamma' = \left( \begin{smallmatrix} A' & B' \\ C & 0
\end{smallmatrix} \right)$, using (\ref{eqn-conditions}) to get $A' = 0;$ hence $B' = -
\tr{C}^{-1}.$  By Lemma \ref{lem-coboundary} it suffices to find $Z \in \mathfrak h_2$ so that
$\widehat{Z} = \gamma' Z.$  Such a point is given by $Z = i\left( \begin{smallmatrix} y_1
 & 0\\ 0 & y_2 \end{smallmatrix} \right)$ where $y_1,y_2>0$ are chosen so that $y_1y_2 =
-d/\lambda^2.$

Now suppose instead that $\lambda = 0.$  Set
\begin{equation*}
\gamma' = \widehat h \gamma h^{-1} = \left( \begin{matrix} A & 0 \\ 0 & \tr{A}^{-1} 
\end{matrix} \right)\text{ where }
h = \left( \begin{matrix} I & \frac{1}{2}A^{-1} B \\ 0 & I \end{matrix} \right).
\end{equation*}
It suffices to find $Z \in \mathfrak h_2$ so that $\gamma' \cdot Z = \widehat{Z}$, that is,
\begin{equation}\label{eqn-needthis}
AZ \tr{A} = \beta \overline{Z} \tr{\beta}^{-1}.
\end{equation}
In fact we will find such an element $Z = iY$ where $Y\in \mathbf{GL}(2,\mathbb R)$ is
positive definite and symmetric.

Since $\beta^2 = dI,$ (\ref{eqn-conditions}) gives $(\beta A)^2 = dI$ hence the matrix
\begin{equation*}
A' = \frac{1}{\sqrt{-d}}\beta A \end{equation*}
has characteristic polynomial $x^2 +1$ so it is $\mathbf{GL}(2,\mathbb R)$-conjugate to the
matrix $j=\left( \begin{smallmatrix} 0 & 1 \\ -1 & 0 \end{smallmatrix} \right),$ say,
$w A' w^{-1} = j.$  Let 
\begin{equation*}
W = \left( \begin{matrix} w & 0 \\ 0 & \tr{w}^{-1} \end{matrix} \right) \text{ and }
J = \left( \begin{matrix} j & 0 \\ 0 & j \end{matrix} \right) \end{equation*}
be the corresponding elements of the symplectic group.  Take $Z = W^{-1} \cdot iI =
iw^{-1}\tr{w}^{-1}.$  Then
\[
A' Z \tr{A}' = w^{-1}jwZ\tr{w}\tr{j}\tr{w}^{-1} = w^{-1}iI \tr{w}^{-1} = Z = - \overline{Z}.\]
Hence
\[ A Z \tr{A} = -d \beta^{-1} (-\overline{Z}) \tr{\beta}^{-1} = \beta \overline{Z}
\tr{\beta}^{-1}.  \qed\]   

\begin{cor}\label{prop-fgamma}
Let $\Gamma \subset \mathbf{Sp}(4,\mathbb R)$ be a torsion-free arithmetic group such that
$\widehat{\Gamma} = \Gamma.$  Fix $\gamma \in \Gamma.$  Then the following statements are
equivalent: \begin{enumerate}
\item $E_{\gamma} \ne \phi$
\item $f_{\gamma}$ is a cocycle (i.e. $\gamma\widehat{\gamma}= I$). 
\item There exists $g \in \mathbf{Sp}(4,\mathbb R)$ such that $\gamma = \widehat{g} g^{-1}.$
\end{enumerate}
In this case, $E_{\gamma} = g\cdot \phi(\mathbf H_3).$  If $\gamma_1,\gamma_2 \in \Gamma$ are
distinct then $E_{\gamma_1} \cap E_{\gamma_2} = \phi.$
\end{cor}

\subsection{Proof}
Suppose $Z \in E_{\gamma}.$  Then $\widehat{\gamma} \gamma Z = \widehat{\gamma Z} = Z$ so
$\widehat{\gamma}\gamma$ is torsion.  By hypothesis this implies $\widehat{\gamma}\gamma = I$
hence (1) implies (2).  If $f_{\gamma}$ is a 1-cocycle then by Proposition \ref{prop-Galois}
it is a coboundary, so (2) implies (3).  If $\gamma = \widehat{g} g^{-1}$ then $E_{\gamma} =
g\cdot \phi(\mathbf H_3)$ so (3) implies (1).  Finally, if $Z \in E_{\gamma_1} \cap
E_{\gamma_2}$ then $\gamma_1 Z = \widehat{Z} = \gamma_2 Z$ so $\gamma_2^{-1} \gamma_1$ fixes
$Z$, hence $\gamma_1 = \gamma_2.$  \qed

\subsection{}  If $\Gamma \subset \mathbf{Sp}(4,\mathbb Z)$ is a subgroup, define
\[ \widetilde{\Gamma} = \left\{ g \in \mathbf{Sp}(4,\mathbb Z)|\
\widehat{g} g^{-1} \in \Gamma \right\}. \]
If $\widehat{\Gamma} = \Gamma$ then $\widehat{\widetilde{\Gamma}} = \widetilde{\Gamma}.$  It
is easy to see that $\widetilde{\Gamma}_N$ is a group and that $\Gamma_N$ is a normal
subgroup.

\begin{prop}  Assume that $\mathcal O_d$ is a principal ideal domain (cf. \S
\ref{subsec-PID}).  Let $N \ge 3.$  If $d \equiv 1 \mod 4,$  assume also that $N$ is even.
Then the short exact sequence
\begin{equation}\label{eqn-shortexact}
1 \longrightarrow \Gamma_N \longrightarrow \widetilde{\Gamma}_N \longrightarrow
\Gamma_N \backslash \widetilde{\Gamma}_N \longrightarrow 1 \end{equation}
induces a bijection
\begin{equation}\label{eqn-doublecoset}
H^1(\mathbb C/ \mathbb R, \Gamma_N) \cong \Gamma_N \backslash \widetilde{\Gamma}_N /
\phi(\mathbf{SL}(2,\mathcal O_d)).  \end{equation}\end{prop}

\subsection{Proof}
  We claim that $H^1(\mathbb C/ \mathbb R, \Gamma_N) \to H^1 (\mathbb C/ \mathbb R,
\widetilde{\Gamma}_N)$ is trivial.  For, if $\gamma \in \Gamma_N$ and if $f_{\gamma}$ is a
1-cocycle then by Corollary \ref{prop-fgamma} the set $E_{\gamma} \ne \phi.$  By Proposition 
\ref{prop-application} (below), there exists $g \in \mathbf{Sp}(4,\mathbb Z)$ so that $\gamma
=\widehat{g}g^{-1}$, which proves the claim.  The long exact cohomology sequence associated to
(\ref{eqn-shortexact}) is
\begin{equation*}\begin{matrix}
H^0(\Gamma_N) &\to& H^0(\widetilde{\Gamma}_N) &\to& H^0(\Gamma_N \backslash
\widetilde{\Gamma}_N) &\to & H^1(\Gamma_N) &\to& 1 \\
|| && || && || \\
\phi(\Lambda_N) && \phi(\Lambda) && \Gamma_N \backslash
\widetilde{\Gamma}_N
\end{matrix}\end{equation*} 
which completes the proof of (\ref{eqn-doublecoset}).  \qed

\section{Anti-holomorphic multiplication}\label{sec-anticomplex}

\subsection{}\label{subsec-PP}
Recall \cite{Lange} that a symplectic form $Q$ on $\mathbb C^2$ is {\it compatible} with the
complex structure if $Q(iu,iv) = Q(u,v)$ for all $u,v\in \mathbb C^2.$  A compatible form $Q$
is {\it positive} if the symmetric form $R(u,v) = Q(iu,v)$ is positive definite.  If $Q$ is
compatible and positive then it is the imaginary part of a unique positive definite Hermitian
form $H = R+iQ.$  Let $L\subset \mathbb C^2$ be a lattice and let $H = R+iQ$ be a positive
definite Hermitian form on $\mathbb C^2.$  A basis of $L$ is {\it symplectic} if the matrix
for $Q$ with respect to this basis is $\left( \begin{smallmatrix} 0 & I \\ -I & 0
\end{smallmatrix} \right).$  The lattice $L$ is symplectic if it admits a symplectic basis.
A principally polarized abelian surface is a pair $(A = \mathbb C^2/L, H = R+iQ)$ where $H$ is
a positive definite Hermitian form on $\mathbb C^2$ and where $L \subset \mathbb C^2$ is a
symplectic lattice relative to $Q = \text{Im}(H).$ 

Each $Z \in \mathfrak {h}_2$ determines a principally polarized abelian surface $(A_Z, H_Z)$
as follows.  Let $Q_0$ be the standard symplectic form on $\mathbb R^4 = \mathbb R^2 \oplus
\mathbb R^2$ with matrix $\left( \begin{smallmatrix} 0 & I \\ -I & 0 \end{smallmatrix}
\right)$ (with respect to the standard basis of $\mathbb R^2 \oplus \mathbb R^2$).  Let $F_Z:
\mathbb R^2 \oplus \mathbb R^2 \to \mathbb C^2$ be the real linear mapping with matrix
$(Z,I)$, that is,
\[ F_Z\left(\begin{smallmatrix} x \\ y \end{smallmatrix} \right ) = Zx + y.\]
Then $(F_Z)_*(Q_0)$ is a compatible, positive symplectic form and $L_Z = F_Z(\mathbb Z^2
\oplus\mathbb Z^2)$ is a symplectic lattice with symplectic basis $F_Z(\text{standard
basis}).$  The Hermitian form corresponding to $Q_Z$ is 
\[H_Z(u,v) = Q_Z(iu,v) + i Q_Z(u,v)= \tr{u} (\text{Im}(Z))^{-1}\bar v \text{ for } u,v \in
\mathbb C^2.\]
The pair $(A_Z = \mathbb C^2/L_Z, H_Z)$ is the desired principally polarized abelian surface.
If $z_1,z_2$ are the standard coordinates on $\mathbb C^2$ then, with respect to the above
symplectic basis of $L$, the differential forms $dz_1, dz_2$ have period matrix $(Z,I).$

The principally polarized abelian surfaces $(A_Z = \mathbb C^2/L_Z, H_Z)$ and $(A_{\Omega} =
\mathbb C^2/L_{\Omega}, H_{\Omega})$ are isomorphic iff there exists a complex linear mapping
$\psi:\mathbb C^2 \to \mathbb C^2$ such that $\psi(L_{\Omega}) = L_Z$ and $\psi_*(H_{\Omega})
=H_Z.$  Set $h = \tr{(F_Z^{-1}\psi F_{\Omega})} = \left( \begin{smallmatrix} A& B\\ C& D
\end{smallmatrix} \right).$  Then: $h\in \mathbf{Sp}(4,\mathbb Z)$, $\Omega = h \cdot Z$, and 
$\psi(M) = \tr{(CZ+D)}M$ for all $M \in \mathbb C^2$,
which is to say that the following diagram commutes:
\begin{equation}\label{diag-abelian1}
\begin{CD}
\left(\begin{smallmatrix} x \\ y \end{smallmatrix} \right)&\quad&
\mathbb R^2 \oplus \mathbb R^2 @>>{F_{h\cdot Z}}> \mathbb C^2 &\quad& M \\
@VVV @VVV @VV{\psi}V @VVV \\
\tr{h}\left(\begin{smallmatrix} x \\ y \end{smallmatrix} \right) &\quad&
\mathbb R^2 \oplus \mathbb R^2 @>>{F_Z}> \mathbb C^2  &\quad& \tr{(CZ+D)}M
\end{CD}\end{equation}
(since $h\cdot Z$ is symmetric). If $Z \in \mathfrak h_2$ and 
$\widehat Z = \beta \overline{Z} \tr{\beta}^{-1}$ then this diagram commutes:
\begin{equation}\label{diag-abelian2}
\begin{CD}
\left(\begin{smallmatrix} x \\ y \end{smallmatrix} \right)&\quad&
\mathbb R^2 \oplus \mathbb R^2 @>>{F_{Z}}> \mathbb C^2 &\quad& M \\
@VVV @VVV @VVV @VVV \\
\tr{N_{\beta}}\left(\begin{smallmatrix} x \\ y \end{smallmatrix} \right) &\quad&
\mathbb R^2 \oplus \mathbb R^2 @>>{F_{\widehat{Z}}}> \mathbb C^2  &\quad& \beta \overline{M}
\end{CD}\end{equation}

\subsection{}\label{subsec-multiplication}  As in \S \ref{subsec-numberfield}, fix a
square-free integer $d < 0$ and let $\mathcal O_d$ denote the ring of
integers in the number field $\mathbb Q[\sqrt{d}].$
Let us say that a principally polarized abelian surface $(A=\mathbb C^2/L, H)$
admits {\it anti-holomorphic multiplication} by the ring $\mathcal O_d$ if there is a
homomorphism $\Psi: \mathcal O_d \to \text{End}_{\mathbb R}(A)$ such that the endomorphism
$\Psi(\sqrt{d})$ is anti-holomorphic and is compatible with the polarization $H.$  To be
precise, suppose $H = R+iQ.$ Then $\Psi$ is an anti-holomorphic multiplication by $\mathcal
O_d$ iff  
\[\kappa = \Psi(\sqrt{d}): \mathbb C^2 \to \mathbb C^2\] 
is complex anti-linear and satisfies \begin{enumerate}
\item $\kappa^2 = dI$ 
\item $Q(\kappa(u),\kappa(v)) = dQ(u,v)$ for all $u,v\in \mathbb C^2$
\item $\kappa(L) \subset L$ and, 
\item if $d \equiv 1 \mod 4$ then $\frac{1}{2}(\kappa +I)(L) \subset L$
\end{enumerate}
in which case, $\Psi$ is determined by $\kappa.$
\begin{lem}\label{lem-ACmult}
Let $\Gamma \subset \mathbf{Sp}(4,\mathbb Z)$ be a torsion-free arithmetic subgroup.   If $d
\equiv 1 \mod 4$ then assume also that $\Gamma
\subset\Gamma(2).$  Suppose $Z \in \mathfrak h_2^{\Gamma}.$  Then there is a unique $\gamma
\in \Gamma$ such that $\widehat{Z} = \gamma \cdot Z.$  Moreover, $\gamma \in \Gamma \cap
\widehat{\Gamma}$, and the mapping
\begin{equation}\label{eqn-kappaZ}
\kappa_Z = F_Z \circ \tr{(}N_{\beta}\gamma) \circ F_Z^{-1}: \mathbb C^2 \to \mathbb C^2
\end{equation}
defines an anti-holomorphic multiplication by $\mathcal O_d$ on the principally polarized
abelian surface $(A_Z,H_Z).$ \end{lem}

\subsection{Proof}  The element $\gamma$ is unique since $\Gamma$ is torsion-free.
Since $\widehat{\gamma}\gamma Z = Z$ we get:  $\widehat{\gamma} \gamma = I.$ Set $\eta =
\tr{(N_{\beta}\gamma)}.$   Then  
\[
\eta^2 = \tr{(N_{\beta} \gamma N_{\beta} \gamma)} = \tr{(N_{\beta} \gamma d N_{\beta}^{-1}
\gamma)} = d\tr{(\widehat \gamma \gamma)} = dI\]
so the same is true of $\kappa_Z.$    Also
\[ \eta = \tr{\gamma} \left( \begin{matrix} \tr{\beta} & 0 \\ 0 & {\beta}^{-1} \end{matrix}
\right) \left( \begin{matrix} I & 0 \\ 0 & dI \end{matrix} \right).\]
The first two factors are in $\mathbf{Sp}(4,\mathbb R)$ so
\[ Q_0 (\eta u, \eta v) = d Q_0(u,v) \]
for all $u,v\in \mathbb R^4.$  Hence $Q_Z(\kappa_Z u, \kappa_Z v) = d Q_Z(u,v)$ for all
$u,v\in\mathbb C^2$ which verifies conditions (1) and (2) above.  Condition (3) holds since
$\eta $ preserves the lattice $\mathbb Z^2 \oplus \mathbb Z^2.$  Now suppose that $d \equiv 1
\mod 4$ and that $\gamma \equiv I \mod 2.$  Then $\tr{(N_{\beta}\gamma)} +I \equiv 0 \mod 2$
(since $\beta \equiv I \mod 2$), which shows that  $\frac{1}{2}(I+\eta)$
preserves the lattice $\mathbb Z^2 \oplus \mathbb Z^2,$ and verifies (4).

Finally we check that $\kappa_Z:\mathbb C^2 \to \mathbb C^2$ is anti-linear.  Let $\gamma =
\left(\begin{smallmatrix} A & B \\ C& D \end{smallmatrix}\right).$ By
(\ref{diag-abelian1}) and (\ref{diag-abelian2}) the following diagram commutes,
\begin{equation}\label{eqn-diagkappa}\begin{CD}
\mathbb R^2\oplus \mathbb R^2 @>>{F_Z}> \mathbb C^2\quad  && M\\
@V{\tr{N_{\beta}}}VV @VVV @VVV\\
\mathbb R^2 \oplus \mathbb R^2 @>>F_{\widehat{Z}}> \mathbb C^2 \quad &&\beta \overline{M} \\
@V{\tr{\gamma}}VV @VVV @VVV\\
\mathbb R^2 \oplus \mathbb R^2 @>>{F_{\gamma^{-1}\cdot\widehat{Z}}}> \mathbb C^2 \quad &&
\tr{(CZ+D)}\beta\overline{M}
\end{CD}\end{equation}
But $Z = \gamma^{-1}\cdot  \widehat{Z}$ so the bottom arrow is also $F_Z$.  Then $\kappa_Z$ is
the composition along the right hand vertical column and it is given by $M \mapsto \tr{(CZ +
D)}\beta \overline{M}$ which is anti-linear.  \qed

\section{The Comessatti lemma}\label{sec-Comessatti}
\subsection{}\label{subsec-PID}  Throughout \S \ref{sec-Comessatti} we fix a square-free
integer $d < 0.$   If $d \equiv 1 \mod 4$ set $m = (d-1)/4.$  Recall \cite{Stark} that
$\mathcal O_d\subset \mathbb Q(\sqrt{d})$ is a principal ideal
domain iff $d \in \{-1,-2,-3,$ $-7,-11,$ $-19,-43,$ $-67,-163\}.$  Define the matrix
$N_{\beta}$ by (\ref{eqn-Nbeta}).

 \begin{lem}\label{lem-Omodule} 
Let $M$ be a free $\mathcal O_d$ module of rank 2.  Let $\Psi(x):M \to M$ denote the action by
$x \in \mathcal O_d.$  Then $M$ has a $\mathbb Z$-basis with respect to which the matrix for
$\Psi(\sqrt{d})$ is $\tr{N_{\beta}}.$
\end{lem}
\subsection{Proof}
Set $M = \mathcal O_d v \oplus \mathcal O_d
w.$  If $ d \not\equiv 1  \mod 4$ then $\mathcal O_d = \mathbb Z[\sqrt{d}]$ and the desired
basis is given by $a_1=v$,$ a_2 = \Psi(\sqrt{d})v$, $b_1 = \Psi(\sqrt{d})w$, and $b_2=w.$

If $d \equiv 1 \mod 4$ let $\alpha = (-1 + \sqrt{d})/2$ and $\alpha' = (1+\sqrt{d})/2.$
Then $\mathcal O_d = \mathbb Z[\alpha] = \mathbb Z[\alpha']$,  $m = \alpha^2+\alpha =
(\alpha')^2 - \alpha'$ and $\sqrt{d} = 2\alpha+1 = 2\alpha'-1.$  The desired basis is given
by $a_1=v$, $a_2 = \Psi(\alpha)v$, $b_1 = \Psi(\alpha')w,$ and $b_2=w.$  \qed

\subsection{}\label{subsec-newperiod}
  Let $A=\mathbb C^2/L$ be a complex torus.  If $\omega_1,\omega_2$ form a
basis for the space of holomorphic 1-forms and if $v_1,v_2,v_3,v_4$ are a basis for $L$, then
the corresponding period matrix $\Omega$ has entries $\Omega_{ij} = \int_{v_j}\omega_i.$  If 
$v'_i = \sum_jA_{ij}v_j$ and if $\omega'_i = \sum_j B_{ij}\omega_j$ then the resulting period
matrix is 
\begin{equation}\label{eqn-newperiod}
\Omega' = B\Omega \tr{A}.\end{equation}
The following proposition is an analog of the lemma (\cite{Silhol2}, \cite{Comessatti}) of
Comessatti and Silhol.

\begin{prop}\label{prop-Comessatti}  Suppose $\mathcal O_d$ is a principal ideal domain.
Let $(A = \mathbb C^2/L, H=R+iQ)$ be a principally polarized abelian surface which admits
an anti-holomorphic multiplication $\kappa: \mathbb C^2 \to \mathbb C^2$ by $\mathcal O_d.$
Then there exists a basis for the holomorphic 1-forms on $A$ and a symplectic basis for the
lattice $L$ such that the resulting period matrix is $ \left( Z, I \right) $
for some $Z \in \phi(\mathbf {H_3}) \subset \mathfrak h_2,$ and such that the matrix for
$\kappa$ with respect to this basis (of $L$) is $\tr{N}_{\beta}.$  
\end{prop}

\subsection{Proof}   The mapping $\kappa$
preserves $L$ and defines on $L$ the structure of an $\mathcal O_d$ module.  Since $L$ is
torsion-free and since $\mathcal O_d$ is principal, it follows that $L$ is a free $\mathcal
O_d$ module of rank 2, to which we may apply Lemma \ref{lem-Omodule}.

First consider the case $ d \not\equiv 1 \mod 4.$  Let $\left\{a_1,a_2,b_1,b_2\right\}$
be the basis of $L$ constructed in Lemma \ref{lem-Omodule}.  Then $\kappa(a_1)=a_2$,
$\kappa(a_2)=da_1$, $\kappa(b_1)=db_2$, and $\kappa(b_2)=b_1.$  We will now show that this
basis is symplectic with respect to the symplectic form $Q$, or else it can be modified to
give a symplectic basis.  Let $x=Q(a_1,b_1)\in \mathbb Z$ and $y=Q(a_1,b_2)\in \mathbb Z.$
Since (\S \ref{subsec-multiplication}) $Q(\kappa u, \kappa v) = dQ(u,v)$, we find that the
matrix for $Q$ is $\left( \begin{smallmatrix} 0 & M \\ -M & 0 \end{smallmatrix} \right)$ where
$M =\left( \begin{smallmatrix} x & y \\ dy & x \end{smallmatrix} \right).$  Since
$\{a_1,a_2,b_1,b_2\}$ form a basis of $L$, we have: $\pm1 = \det(M)=x^2 - dy^2 >0.$ Hence this
determinant is 1 and there are four possibilities \begin{itemize}
\item $x=1$ and $y=0$
\item $x=-1$ and $y=0$
\item $x=0$ and $y=1$ and $d=-1$
\item $x=0$ and $y=-1$ and $d=-1.$
\end{itemize} 
In the second case, replace $a_1$ by $a'_1 = -a_1$ and replace $a_2$ by $a'_2 = -a_2.$  In the
third case replace $a_1$ by $a'_1 = a_2$ and replace $a_2$ by $a'_2 = -a_1.$  In the fourth
case replace $a_1$ by $a'_1 = -a_2$ and replace $a_2$ by $a'_2 = a_1.$  Then in all four
cases, the resulting basis $\{a'_1,a'_2,b'_1=b_1, b'_2=b_2\}$ is symplectic and the matrix for
$\kappa$ is $\tr{N_{\beta}}.$

Choose a holomorphic 1-form $\omega_1 \in \Gamma(\Omega^1_A)$ such
that $\int_{b'_1}\omega_1 = 1$  and  $\int_{b'_2}\omega_1 = 0.$  Define $\omega_2 =
\overline{\kappa^*\omega_1}$ to be the complex conjugate of the pullback of $\omega_1$ under
the mapping $\kappa.$  Set $z= \int_{a'_1}\omega_1$ and $r = \int_{a'_2}\omega_1.$  Then
\begin{align*}
\int_{b'_i}\omega_2 = \int_{\kappa(b'_i)} \overline{\omega_1} = 
&\begin{cases} 0 &\text{ if }i=1 \\ 1 &\text{ if } i=2 \end{cases}\\
 \int_{a'_i} \omega_2 = \int_{\kappa(a'_i)} \overline{\omega_1} =
&\begin{cases} \bar r &\text{ if } i=1 \\ d\bar z &\text{ if } i=2 \end{cases} \end{align*}
Therefore the period matrix for $\{\omega_1,\omega_2,a'_1,a'_2,b'_1,b'_2\}$ is $\left( \Omega,
I\right)$ where $\Omega = \left( \begin{smallmatrix} z & r \\ \bar r & d\bar z
\end{smallmatrix} \right).$  The Riemann relations imply that this matrix is symmetric (so
$\bar r = r$) and that $\text{Im}(\Omega) > 0.$  This completes the proof if $d\not\equiv 1\
\mod 4.$  

Now consider the case $ d \equiv 1 \mod 4.$  Let $\left\{ a_1,a_2,b_1,b_2\right\}$ be the
basis of $L$ which was constructed in Lemma \ref{lem-Omodule}.  The matrix for $\kappa$ with
respect to this basis is $\tr{N}_{\beta}.$  Set
\begin{equation}\label{eqn-basischange}\begin{array}{rrrrrr}
a'_1 &=&  \frac{1}{2}a_1\\
a'_2 &=& \frac{1}{2}a_1 &+& a_2\\
b'_1 &=& 2b_1 &-& b_2\\
b'_2 &=&  && b_2
\end{array}\end{equation}
 Then 
\begin{equation}\label{eqn-kappa}
\kappa(a'_1)=a'_2,\ \kappa(a'_2)=da'_1,\  \kappa(b'_1)=db'_2,\text{ and }
\kappa(b'_2)=b'_1.\end{equation}  
We will modify this basis so as to obtain a symplectic basis of $L$.  Set $x' = Q(a'_1,b'_1)
\in \frac{1}{2}\mathbb Z$ and $y' = Q(a'_1,b'_2) \in \frac{1}{2}\mathbb Z.$  The matrix for
$Q$ with respect to the basis $\{ a'_1,a'_2,b'_1,b'_2\}$ of $\mathbb C^2$ is $\left(
\begin{smallmatrix} 0 & M' \\ -M' & 0 \end{smallmatrix} \right)$ where $M' = \left(
\begin{smallmatrix} x' & y' \\ dy' & x' \end{smallmatrix} \right).$  The 
determinant of (\ref{eqn-basischange}) is 1 so $\pm 1 = \det(M') = (x')^2 - d(y')^2 >0.$
Since $d \le -3$ there are the following possibilities, \begin{itemize}
\item $x'=1$ and $y'= 0$
\item $x'=-1$ and $y'=0$
\item $x' = \pm \frac{1}{2}$, $y' = \pm \frac{1}{2}$, and $d = -3.$
\end{itemize}
In the second case, replace $a'_1$ by $-a'_1$ and $a'_2$ by $-a'_2.$  In the third case
replace $a'_1$ by $ua'_1$ and $a'_2$ by $ua'_2$ where 
\[ u = \pm \left( \frac{1 \pm \sqrt{-3}}{2} \right) \]
is the appropriately chosen unit in $\mathcal O_{-3}.$  In all three cases the resulting
basis, which we will continue to denote by $\{ a'_1,a'_2,b'_1,b'_2\}$ is a symplectic basis
for $\mathbb C^2$ such that (\ref{eqn-kappa}) holds, although it is not a basis of the lattice
$L$.

Choose a holomorphic 1-form $\omega'_1$ so that $\int_{b'_1}\omega'_1=1$ and
$\int_{b'_2}\omega'_1=0.$  Set $\omega'_2 =\overline{\kappa^*\omega'_1}.$   Set $z =
\int_{a'_1}\omega'_1$ and $r = \int_{a'_2}\omega'_1.$  Then the period matrix for
$\{\omega'_1,\omega'_2,a'_1,a'_2,b'_1,b'_2\}$ is $\left(\Omega,I \right)$ where $\Omega =
\left(\begin{smallmatrix} z & \bar r \\ r & d\bar z \end{smallmatrix} \right)$ as
in the case $d \not\equiv 1 \mod 4.$  Now define a further change of coordinates,
\begin{equation*}
\begin{array}{rrrrrrrrr}
\omega^{\prime\prime}_1 &=& 2\omega_1 & \\
\omega^{\prime\prime}_2 &=& -\omega_1 &+&\omega_2\\
a^{\prime\prime}_1 &=&2a'_1 && &=& ua_1\\
a^{\prime\prime}_2 &=& -a'_1 &+&a'_2 &=& ua_2 \\
b^{\prime\prime}_1 &=& \frac{1}{2} b'_1 &+&\frac{1}{2} b'_2 &=& b_1 \\
b^{\prime\prime}_2 &=& &&b'_2 &=& b_2 
\end{array}\end{equation*} 
(Here, $u$ is a unit in $\mathcal O_d$ and it equals 1 except possibly in the case $d=-3.$)
Then $\{a_1^{\prime\prime}, a_2^{\prime\prime}, b_1^{\prime\prime}, b_2^{\prime\prime} \}$ is
an integral symplectic basis for $L.$  By (\ref{eqn-newperiod}) the period matrix for $\{
\omega_1^{\prime\prime},$ $\omega_2^{\prime\prime},$ $a_1^{\prime\prime},$
$a_2^{\prime\prime},$ $b_1^{\prime\prime},$ $b_2^{\prime\prime}\}$ is given by
\begin{equation*}
\sigma \left( \Omega, I \right) \left( \begin{matrix} \tr{\sigma} & 0 \\ 0 & {\sigma}^{-1} 
\end{matrix} \right) =\left(\sigma\Omega\tr{\sigma}, I \right)\end{equation*}
where $\sigma = \left( \begin{smallmatrix}2&0\\-1&1\end{smallmatrix} \right).$
It follows that $\Omega$ is symmetric, $\text{Im}(\Omega) >0,$ and the matrix for $\kappa$
with respect to the basis $\left\{ a_1^{\prime\prime}, a_2^{\prime\prime}, b_1^{\prime\prime},
b_2^{\prime\prime} \right\}$ is $\tr{N}_{\beta}.$  \qed

\subsection{}
The main application of the Comessatti lemma is to strengthen the conclusion of Proposition
\ref{prop-fgamma}.  Fix a square-free integer $d < 0$ and let $\mathcal O_d$ be the ring of
integers in $\mathbb Q(\sqrt{d}).$
Let $\Gamma \subset \mathbf{Sp}(4,\mathbb Z)$ be a torsion free subgroup such that
$\widehat\Gamma = \Gamma.$  If $d \equiv 1 \mod 4$ then suppose also that $\Gamma \subset
\Gamma(2).$
\begin{prop}\label{prop-application}
Let $\gamma \in \Gamma$ and suppose $\mathcal O_d$ is a principal ideal domain.  Then
$E_{\gamma} \ne \phi$ if
and only if there exists $g \in \mathbf{Sp}(4,\mathbb Z)$ such that $\gamma = \widehat g
g^{-1}.$ In this case the element $g$ is uniquely determined up to multiplication from the
right by an element of $\phi(\mathbf{SL}(2,\mathcal O_d)),$ and
\[ E_{\gamma} = g \cdot \phi(\mathbf H_3).\]
\end{prop}

\subsection{Proof}
Suppose $\gamma = \widehat g g^{-1}$ for some $g\in \mathbf{Sp}(4,\mathbb Z).$  Then $Z \in
E_{\gamma}$ iff $\gamma \cdot Z = \widehat Z$ iff $\widehat g g^{-1}\cdot Z = \widehat Z$ iff
$g^{-1}Z \in \phi(\mathbf H_3)$ since it is fixed under the $\widehat{\ } $ involution.  This
proves the ``if'' part of the first statement.  

Now suppose that $E_{\gamma} \ne \phi.$  Choose $Z \in E_{\gamma}$ which is not fixed by any
element of $\mathbf{Sp}(4,\mathbb Z)$ other than $\pm I.$  Such points exist and are even
dense in $E_{\gamma}$ since, by Proposition \ref{prop-fgamma} the set $E_{\gamma}$ is a
translate of $\phi(\mathbf H_3).$  Let $L_Z= F_Z(\mathbf Z^2 \oplus \mathbf Z^2)$ be the
corresponding lattice in $\mathbb C^2$ and let $(A_Z = \mathbb C^2/L_Z, H_Z)$ be the
corresponding polarized abelian surface.  By Lemma \ref{lem-ACmult} this variety admits
anti-holomorphic multiplication $\kappa_Z$ by the ring $\mathcal O_d.$  By Proposition
\ref{prop-Comessatti} there is a basis for the holomorphic 1-forms on $A_Z$ and a symplectic
basis for $L_Z$ so that the resulting period matrix is $\left(\Omega,I \right)$ where
\[ \Omega= \sigma \left( \begin{matrix} z & r \\ r & d\bar z \end{matrix}
\right)\tr{\sigma}\ \text{  with  } r \in \mathbb R \text{  and  }
\text{Im}\left( \begin{matrix} z & r \\ r & d\bar z \end{matrix}
\right)>0.\]
Then $(\mathbb C^2/L_Z, H_Z)$ and $(\mathbb C^2/L_{\Omega}, H_{\Omega})$ are isomorphic
principally polarized abelian surfaces, so by \S \ref{subsec-PP} there exists a linear
mapping $\psi:\mathbb C^2 \to \mathbb C^2$ so that $\psi(L_{\Omega}) = L_Z$ and
$\psi_*(H_{\Omega}) = H_Z.$  Define $h = \tr{(F_Z^{-1}\psi F_{\Omega})}
 \in \mathbf{Sp}(4,\mathbb Z).$  Then $\Omega = h \cdot Z$.  By Proposition
\ref{prop-embeddings}, $\Omega =
\widehat{\Omega} = \widehat h \widehat{Z} = \widehat h \gamma Z.$ So $h^{-1} \widehat h \gamma
= \pm I.$  If the plus sign occurs, then take $g = h^{-1}$ to get $\gamma = \widehat g
g^{-1}.$  If the minus sign occurs, take $g = h^{-1}\left(\begin{smallmatrix} 0 & I \\ -I & 0
\end{smallmatrix} \right)$ to get $\gamma = \widehat g g^{-1}.$

Finally suppose $\gamma = \widehat{h} h^{-1}$ for some $h \in \mathbf{Sp}(4,\mathbb Z).$  Then
the element $ a = h^{-1} g \in \mathbf{Sp}(4,\mathbb Z)$ satisfies $\widehat{a} = a.$  By
Propositions \ref{prop-embeddings} and \ref{prop-arithmeticgroup},
\[ a \in \mathbf{Sp}(4,\mathbb Z) \cap \phi(\mathbf{SL}(2,\mathbb C)) =
\phi(\mathbf{SL}(2,\mathcal O_d)).\] 
\qed.

\section{Arithmetic quotients}
As in \S \ref{subsec-numberfield} fix a square-free integer $d < 0,$ let $\mathcal O_d$ be the
ring of integers in $\mathbb Q(\sqrt{d}).$  To simplify notation  we will write $\Lambda =
\mathbf{SL}(2,\mathcal O_d).$  As in \S \ref{subsec-arithmeticsubgroups}, fix an integer $N
\ge 1,$ let $\Gamma_N = \Gamma(N) \cap \widehat{\Gamma}(N)$ and denote by $\Lambda_N$ the
principal level subgroup
\[ \Lambda_N = \mathbf{SL}(2,\mathcal O_d)(N).\]
As in Proposition \ref{prop-embeddings}, the group $\mathbf{SL}(2,\mathbb C)$ acts on
hyperbolic space $\mathbf H_3$ in a way which is compatible with the embeddings $\phi:
\mathbf{SL}(2,\mathbb C) \to \mathbf{Sp}(4,\mathbb R)$ and $\phi: \mathbf H_3 \to \mathfrak
h_2.$  Let $W = \Lambda_N \backslash \mathbf H_3$ and let $X = \Gamma_N \backslash \mathfrak
h_2.$  Let $\pi: \mathbf H_3 \to W$ and $\pi: \mathfrak h_2 \to X$ denote the projections.

For any $g \in \widetilde{\Gamma}_N$ define
\[ {}^g\phi: W = \Lambda_N \backslash \mathbf H_3 \to \Gamma_N \backslash \mathfrak h_2 = X\]
by $\Lambda_N w \mapsto \Gamma_N g \phi(w).$

\begin{lem}\label{lem-injective}  Let $g \in \widetilde{\Gamma}_N$ and set $\gamma =
\widehat{g}g^{-1}\in
\Gamma_N.$  Then the mapping ${}^g\phi$ is well defined and injective. Its image
\[  {}^g\phi(W) = \pi(E_{\gamma})\]
coincides with the projection of $E_{\gamma}.$  If $h,g \in \widetilde{\Gamma}_N$ determine
the same double coset 
\[ \Gamma_N h \phi(\Lambda) = \Gamma_N g \phi(\Lambda) \in \Gamma_N \backslash
\widetilde{\Gamma}_N / \phi(\Lambda)\] 
then ${}^g\phi(W) = {}^h\phi(W).$
\end{lem}

\subsection{Proof}
The mapping ${}^g\phi$ is well defined, for if $\lambda \in \Lambda_N$ then
\[ {}^g\phi(\lambda w) = \Gamma_N g \phi(\lambda)g^{-1} g\phi(w) = \Gamma_N g \phi(w) 
= {}^g\phi(w)\]
because $\phi(\lambda) \in \Gamma_N$ which is normal in $\widetilde{\Gamma}_N.$  It is easy to
see that $E_{\gamma} = g \cdot \phi(\mathbf H_3) \subset \mathfrak h_2$, so we have a
commutative diagram
\begin{equation*}\begin{CD}
\mathbf H_3 @>{g\cdot \phi}>> E_{\gamma} & \subset\ & \mathfrak h_2 \\
@V{\pi}VV @VVV @VV{\pi}V \\
W @>{{}^g\phi}>> \pi(E_{\gamma})\ & \subset\ & X \end{CD}  \end{equation*}
from which it follows that ${}^g\phi(W) = \pi(E_{\gamma}).$  We need to show that ${}^g\phi$
is injective.

We claim the normalizer $N_{\gamma}$ of $E_{\gamma}$ in $\Gamma_N$ is $g \phi(\Lambda_N)
g^{-1}.$   For, if $h \in \Gamma_N$, $x \in E_{\gamma}$, and $hx\in E_{\gamma}$
then $\widehat h \widehat x = \gamma hx$ so $h^{-1} \gamma^{-1} \widehat h \gamma \in
\Gamma_N$ fixes $x$.  This holds for all $x \in E_{\gamma}$ so
 $\widehat h \gamma = \gamma h$, or $\widehat g^{-1} \widehat h \widehat g = g^{-1} h g.$
This implies that $g^{-1} h g \in \phi(\mathbf{SL}(2,\mathbb C))$ hence $h \in g
\phi(\mathbf{SL}(2,\mathbb C))g^{-1} \cap \Gamma_N$ which proves the claim.

Now let $w_1,w_2 \in \mathbf H_3$ and suppose that ${}^g\phi(\Lambda_Nw_1) = {}^g \phi(
\Lambda_Nw_2)\in \pi(E_{\gamma}).$ Since $\pi(E_{\gamma}) = N_{\gamma} \backslash E_{\gamma}$
there exists $h \in N_{\gamma}$ so that $h g \phi(w_1) = g \phi(w_2).$  By the preceding
claim, $h = g \phi(\lambda) g^{-1}$ for some $\lambda \in \Lambda_N,$ which implies $g
\phi(\lambda w_1) = g \phi(w_2).$  But $\phi$ is injective so $w_2 = \lambda w_1$ and
$\Lambda_N w_2 = \Lambda_N w_1 \in W.$  

Finally, if $h = \gamma g \phi(\lambda)$ represents the same double coset as $g$ (where
$\gamma \in \Gamma_N$ and $\lambda \in \Lambda$) then \vskip-25pt 
\[ {}^h\phi(\Lambda_Nw) = \Gamma_N \gamma g \phi(\lambda) \phi(w) = {}^g\phi(\lambda w) \]
which shows that ${}^h\phi(W) = {}^g\phi(W).$   \qed   

The involution $\ \widehat{}\ : \mathfrak h_2 \to \mathfrak h_2$ passes to an anti-holomorphic
involution on $X$ and hence defines a real structure on $X.$

\begin{thm}\label{thm-realpoints}
Suppose $\mathcal O_d$ is a principal ideal domain.  Fix $N \ge 3.$  If $d \equiv 1 \mod
4$ assume also that $N$ is even.  Then the set of real points $X_{\mathbb R}$ is the disjoint
union
\begin{equation}\label{eqn-realpoints}
X_{\mathbb R} = \coprod_g {}^g\phi(W) \end{equation}
of finitely many copies of $W = \Lambda_N \backslash \mathbf H_3,$ indexed by elements
\[ g \in \Gamma_N \backslash \widetilde{\Gamma}_N / \phi(\Lambda) = H^1(\mathbb C/\mathbb R,
\Gamma_N).\]
The copy indexed by $g \in \widetilde{\Gamma}$ corresponds to the cohomology class
$f_{\gamma}$ where $\gamma = \widehat{g} g^{-1} \in \Gamma_N.$
\end{thm}

\subsection{Proof}
A point $x \in X$ is real iff it is the image of a $\Gamma_N$-real point $Z \in \mathfrak
h_2.$  So 
\[ X_{\mathbb R} = \bigcup_{\gamma \in \Gamma_N} \pi(E_{\gamma}).\]
By Lemma \ref{lem-injective}, for each $g \in \widetilde{\Gamma}_N$ the mapping ${}^g\phi:W
\to X_{\mathbb R}$ is injective.  By Proposition \ref{prop-application} the union
(\ref{eqn-realpoints}) contains $X_{\mathbb R}.$  It remains to show that the individual
components are disjoint.

Let $g_1,g_2 \in \widetilde{\Gamma}_N$ and set $\gamma_i = \widehat{g}_i g_i^{-1} \in
\Gamma_N.$  Suppose the corresponding components have a nontrivial intersection, say
\[ x' \in \pi(E_{\gamma_1}) \cap \pi(E_{\gamma_2}).\]
We claim that there exists $\gamma \in \Gamma_N$ so that 
\begin{equation}\label{eqn-sameimage}
\widehat{\gamma} \gamma_1 = \gamma_2 \gamma \end{equation}
and hence $\pi(E_{\gamma_1}) = \pi(E_{\gamma_2}).$  There is a lift $x \in E_{\gamma_1}$ of
$x'$ so that $\gamma x \in E_{\gamma_2}$ for some $\gamma \in \Gamma_N.$  This implies
$\widehat \gamma \gamma_1 x = \widehat \gamma \widehat x = \gamma _2 \gamma x.$  Therefore
$\gamma^{-1} \gamma_2^{-1} \widehat \gamma \gamma_1\in \Gamma_N$ fixes $x.$  But $\Gamma_N$ is
torsion-free, so $\widehat \gamma \gamma_1 = \gamma_2 \gamma$  as claimed.

Equation (\ref{eqn-sameimage}) gives $\widehat{g}_2^{-1} \widehat{\gamma} \widehat{g}_1 =
g_2^{-1} \gamma g_1\in \mathbf{Sp}(4,\mathbb Z)$ which implies by Propositions
\ref{prop-embeddings} and \ref{prop-arithmeticgroup} that $g_2^{-1} \gamma g_1 \in
\phi(SL_2(\mathcal O_d)).$  In other
words, $g_2 = \gamma g_1 \phi(h)$ for some $h \in \mathbf{SL}(2,\mathcal O_d)$, so the
elements $g_1$ and $g_2$ represent the same double coset in $\Gamma_N \backslash
\widetilde{\Gamma} / \phi(\Lambda).$  \qed

\section{The moduli space of abelian surfaces with anti-holomorphic multiplication}
\label{sec-modulispace}
\subsection{Level structures}  Let $(A = \mathbb C^2/L, H=R+iQ)$ be a principally polarized
abelian surface.  A level $N$ structure on $A$ is a choice of basis $\left\{
U_1,U_2,V_1,V_2\right\}$ for the $N$-torsion points of $A$ which is symplectic, in the sense
that there exists a symplectic basis $\left\{ u_1,u_2,v_1,v_2\right\}$ for $L$ such that
\[
U_i \equiv \frac{u_i}{N} \text{ and } V_i \equiv \frac{v_i}{N}\ \mod L\]
(for $i = 1,2$).  For a given leven $N$ structure, such a choice $\left\{
u_1,u_2,v_1,v_2\right\}$ determines a mapping
\begin{equation}\label{eqn-liftmapping}
 F: \mathbb R^2 \oplus \mathbb R^2 \to \mathbb C^2 \end{equation}
such that $F(\mathbb Z^2 \oplus \mathbb Z^2) = L,$ by $F(e_i) = u_i$ and $F(f_i) = v_i$ where
$\left\{ e_1,e_2,f_1,f_2\right\}$ is the standard basis of $\mathbb R^2 \oplus \mathbb R^2.$
The choice $\left\{u_1,u_2,v_2,v_2\right\}$ (or equivalently, the mapping $F$) will be
referred to as a {\it lift} of the level $N$ structure.  It is well defined modulo the
principal congruence subgroup $\Gamma(N)$, that is, if $F': \mathbb R^2 \oplus \mathbb R^2 \to
\mathbb C^2$ is another lift of the level structure, then $F' \circ F^{-1} \in \Gamma(N).$

Fix a square-free integer $d < 0$ and let $(A, H, \kappa)$ be a principally polarized abelian
surface with anti-holomorphic multiplication by $\mathcal O_d$ as in \S
\ref{subsec-multiplication}.  A level $N$ structure $\left\{ U_1,U_2,V_1,V_2 \right\}$ on $A$
is {\it compatible} with  $\kappa$ if for some (and hence for any) lift $F$ of the level
structure, the following diagram commutes $\mod L:$
\begin{equation}\label{diag-level}\begin{CD}
\frac{1}{N}\left( \mathbb Z^2 \oplus \mathbb Z^2\right) @>>F> \frac{1}{N} L\\
@V{\textstyle{\tr{N}_{\beta}}}VV @VV{\textstyle{\kappa}}V \\
\frac{1}{N} \left(\mathbb Z^2 \oplus \mathbb Z^2 \right)@>>{F}> \frac{1}{N} L
\end{CD}\end{equation}

\quash{

\setlength{\arraycolsep}{3pt}
\begin{equation*}
\left.\begin{array}{ccrrrrrrrrrrr}
\k{\ua} &\equiv&  & &\ub  &\quad& \k{\va} &\equiv&  & &d\vb  
\\ \\
\k{\ub} &\equiv& d\ua & &  &\quad& \k{\vb} &\equiv& \va & & 
\end{array}\right\}\mod L
\end{equation*}
if $d \not\equiv 1 \mod 4$ and 
\begin{equation*}\left.\begin{array}{ccrrrrrrrrrrr}
\k{\ua} &\equiv& \ua &+ &2\ub  &\quad& \k{\va} &\equiv& \va &+ &2m\vb 
\\ \\
\k{\ub} &\equiv& 2m\ua &- &\ub  &\quad& \k{\vb} &\equiv& 2\va &- &\vb  
\end{array}\right\} \mod L\end{equation*}
if $d \equiv 1 \mod 4.$

} %% end quash

We will refer to the collection
$ \mathcal A = (A = \mathbb C^2/L, H = R+iQ, \kappa, \left\{ U_i, V_j \right\} )$ as a
principally polarized abelian surface with anti-holomorphic multiplication and 
level $N$ structure.  If $\mathcal A' = (A'=\mathbb C^2/L', H' = R+iQ, \kappa', \left\{
U'_i, V'_j\right\})$ is another such, then an isomorphism $\mathcal A \cong \mathcal A'$ is a
complex linear mapping $\psi: \mathbb C^2 \to \mathbb C^2$ such that $\psi(L)=L',$ $\psi_*(H)
= H',$ $\psi_*(\kappa) = \kappa',$ and such that for some (and hence for any) lift $\left\{
u_1,u_2,v_1,v_2\right\}$ and $\left\{ u'_1,u'_2,v'_1,v'_2\right\}$ of the level structures,
\begin{equation*} 
\psi\left(\frac{\textstyle u_i}{\textstyle N}\right) \equiv \frac{\textstyle u'_i}{\textstyle
N}\ \text{ and }\
\psi\left(\frac{\textstyle v_j}{\textstyle N}\right) \equiv \frac{\textstyle v'_j}{\textstyle
N} \quad
\mod L.
\end{equation*}

\subsection{}  If $Z \in \mathfrak h_2,$ then for any $N \ge 1$ we define the {\it standard
level $N$ structure} on the abelian surface $(A_Z, H_Z)$ to be the basis $\left\{
F_Z(e_i/N), F_Z(f_j/N) \right\} \ \mod L$ where $\left\{
e_1,e_2,f_1,f_2\right\}$ is the standard basis of $\mathbb R^2 \oplus \mathbb R^2.$

\begin{lem}\label{lem-compatibility}
Suppose $Z \in \mathfrak h_2,$  $\gamma \in \mathbf{Sp}(4,\mathbb Z),$ and $\widehat{Z} =
\gamma\cdot Z.$  Let $N \ge 3.$  Then the standard level $N$ structure on
the abelian surface $(A_Z,H_Z)$ is compatible with the anti-holomorphic multiplication
$\kappa_Z$ iff $\gamma \in \Gamma_N.$ 
\end{lem}
\subsection{Proof}  It follows immediately from diagram (\ref{eqn-diagkappa}) that $\gamma \in
\Gamma(N)$ iff the standard level $N$ structure on $(A_Z,H_Z)$ is compatible with $\kappa_Z.$  
Since $\Gamma(N)$ is torsion-free, $\widehat{\gamma} \gamma =I$ which implies $\gamma \in
\widehat\Gamma(N);$ hence $\gamma \in \Gamma_N.$  \qed

  By Lemma \ref{lem-compatibility}, each point $Z \in \mathfrak h_2^{\Gamma(N)}$
determines a principally polarized abelian surface $\mathcal A_Z = \left( A_Z, H_Z,
\kappa_Z, \{ F_Z(e_i/N), F_Z(f_j/N)\} \right)$ with anti-holomorphic
multiplication and (compatible) level $N$ structure.

\begin{thm}\label{thm-main}
 Suppose $\mathcal O_d$ is principal.  Fix $N \ge 3.$  If $d \equiv 1 \mod 4,$
assume also that $N$ is even.   Then the association $Z \mapsto \mathcal A_Z$ determines a one
to one correspondence between the real points (\ref{eqn-realpoints}) $X_{\mathbb R}$ of $X
=\Gamma_N \backslash\mathfrak h_2$ and the moduli space $V(d,N)$ consisting of isomorphism
classes of principally polarized abelian surfaces with anti-holomorphic multiplication by
$\mathcal O_d$ and (compatible) level $N$ structure. \end{thm}  

\subsection{Proof}  A point $x \in X$ is real iff it is the image of a $\Gamma_N$-real point
$Z \in \mathfrak h_2.$  If two $\Gamma_N$-real points $Z, \Omega \in \mathfrak h_2$ determine
isomorphic varieties, say
$\psi:\mathcal A_{\Omega} \cong \mathcal A_Z$ then by (\ref{diag-abelian1}) there exists $h
\in \mathbf{Sp}(4,\mathbb Z)$ such that $\Omega = h\cdot Z.$  Since the isomorphism $\psi$
preserves the level $N$ structures, it follows also from (\ref{diag-abelian1}) that $h \in
\Gamma(N).$  We claim that $h \in \Gamma_N.$  Let $\widehat{Z} = \gamma_Z \cdot Z$ and
$\widehat \Omega = \gamma_{\Omega} \cdot \Omega,$ with $\gamma_Z, \gamma_{\Omega} \in
\Gamma_N.$   Putting diagram (\ref{eqn-diagkappa}) for
$Z$ together with the analogous diagram for $\Omega$ and diagram (\ref{diag-abelian1}), and
using the fact that $\psi_*(\kappa_{\Omega}) = \kappa_Z$ gives a diagram
\begin{equation*} \begin{CD}
\mathbb R^2 \oplus \mathbb R^2 @>>{F_{\Omega}}> \mathbb C^2 \\
@V{\textstyle{\tr{(}\gamma_{\Omega}^{-1}\widehat{h} \gamma_Z)}}VV
@VV{\textstyle{\kappa_Z \psi \kappa_{\Omega}^{-1} = \psi}}V \\
\mathbb R^2 \oplus \mathbb R^2 @>>{F_Z}> \mathbb C^2
\end{CD}\end{equation*}
from which it follows that $\tr{(}\gamma_{\Omega} \widehat{h} \gamma_Z) \in \Gamma(N),$ hence 
$\widehat h \in \Gamma(N)$, hence $h \in \Gamma_N.$

So it remains to show that every principally polarized abelian surface with
anti-hol\-omor\-phic multiplication and level $N$ structure, $\mathcal A = (A, H, \kappa,
\left\{U_i,V_j\right\})$ is isomorphic to some $\mathcal A_Z.$  By the Comessatti lemma
(Proposition \ref{prop-Comessatti}) there exists $Z' \in \mathfrak h_2,$ such that
$\widehat{Z}' = Z',$ (cf. Proposition \ref{prop-embeddings}) and there exists an isomorphism 
\[ \psi': (A_{Z'},H_{Z'},\kappa_{Z'}) \cong (A, H, \kappa)\]
between the principally polarized abelian surfaces with anti-holomorphic multiplication. 
However the isomorphism $\psi'$ must be modified because it does not necessarily take the 
standard level $N$ structure on $(A_{Z'},H_{Z'},\kappa_{Z'})$ to the
given level $N$ structure on $(A, H, \kappa).$  

Choose a lift $\left\{ u_1,u_2,v_1,v_2\right\}$ of the level $N$ structure and let
$F: \mathbb R^2 \oplus \mathbb R^2 \to \mathbb C^2$ be the corresponding mapping
(\ref{eqn-liftmapping}).  Define
\begin{align}
\tr{g}^{-1} &= F^{-1} \circ \psi' \circ F_{Z'} \in \mathbf{Sp}(4,\mathbb Z) \\
Z &= g \cdot Z' \\  \label{eqn-gammaprime}
\gamma &= \widehat{g} g^{-1} = N_{\beta}^{-1} g N_{\beta} g^{-1}.
\end{align}
If $ g = \left(
\begin{smallmatrix} A & B \\ C & D \end{smallmatrix} \right)$ define $\xi: \mathbb C^2 \to
\mathbb C^2$ by $\xi (M) = \tr{(}CZ+D)M.$  Define $\psi = \psi' \circ \xi.$  We will show that
$\gamma \in \Gamma_N,$ that $\widehat{Z} = \gamma \cdot Z,$ and that $\psi$ induces an
isomorphism $\psi: \mathcal A_Z \to \mathcal A$ of principally polarized abelian surfaces with
anti-holomorphic multiplication and compatible level $N$ structures.

In the following diagram, $F$ is the mapping (\ref{eqn-liftmapping}) associated to the lift of
the level $N$ structure.  The bottom square commutes by the definition of $g,$  while the top
square commutes by (\ref{diag-abelian1}).
\begin{equation}\label{diag-abelian1prime}
\begin{CD}
\mathbb R^2 \oplus \mathbb R^2 @>{F_{Z}}>> \mathbb C^2\\
@V{\tr{g}}VV @VV{\xi}V \\
\mathbb R^2 \oplus \mathbb R^2 @>>{F_{Z'}}> \mathbb C^2 \\
@V{\tr{g}^{-1}}VV @VV{\psi'}V \\
\mathbb R^2 \oplus \mathbb R^2 @>{F}>> \mathbb C^2
\end{CD}.\end{equation}

First let us verify that $\xi: (A_Z, H_Z, \kappa_Z) \to (A_{Z'}, H_{Z'}, \kappa_{Z'})$ is an
isomorphism of principally polarized varieties with anti-holomorphic multiplication by
$\mathcal O_d.$  It follows from (\ref{diag-abelian1prime}) that $\xi_*(L_Z) = L_{Z'}$ and
$\xi_*(H_Z)=H_{Z'}.$  We claim that $\xi_*(\kappa_Z) =
\kappa_{Z'},$ that is, $\kappa_{Z'} = \xi \kappa_Z \xi^{-1}.$  But this follows from direct
calculation using $\xi = F_{Z'} \tr{g} F_Z,$ $\kappa_Z = F_Z \tr{(}N_{\beta} \gamma)
F_Z^{-1},$ $\kappa_{Z'} = F_{Z'}\tr{N}_{\beta} F_{Z'}$ and (\ref{eqn-gammaprime}) (and it is
equivalent to the statement that the pushforward by ${} \tr{g}$ of the involution ${}
\tr{(}N_{\beta} \gamma)$ on $\mathbb R^2 \oplus \mathbb R^2$ is the involution ${}
\tr{N}_{\beta}$). It follows that 
\begin{equation}\label{eqn-psikappa}
\psi_*(\kappa_Z) = \kappa.\end{equation}

We claim that the standard level $N$ structure on $(A_Z,H_Z)$ is compatible with $\kappa_Z.$
By construction, the mapping $\psi$ takes the standard level $N$ structure on $(A_Z, H_Z)$ to
the given leven $N$ structure on $(A,H).$  By assumption, the diagram (\ref{diag-level})
commutes $ \mod L.$  By (\ref{diag-abelian1prime}), $F = \psi \circ F_Z.$  Using
(\ref{eqn-psikappa}) it follows that the diagram
\begin{equation*}\begin{CD}
\frac{1}{N} \left(\mathbb Z^2 \oplus \mathbb Z^2\right) @>>{F_Z}> \frac{1}{N} L_Z \\
@V{\tr{N}_{\beta}}VV @VV{\kappa_Z}V \\
\frac{1}{N}\left( \mathbb Z^2 \oplus \mathbb Z^2\right) @>>{F_Z}> \frac{1}{N}L_Z
\end{CD} \end{equation*}
commutes $ \mod {L_Z},$ which proves the claim.  It also follows from Lemma
\ref{lem-compatibility} that $\gamma \in \Gamma_N.$ 

In summary, we have shown that $\left(A_Z,H_Z,\kappa_Z, \left\{ F_Z({e_i}/{N}),
F_Z({f_j}/{N})\right\}\right)$ is a real principally polarized abelian surface with
anti-holomorphic multiplication and (compatible) level $N$ structure, and that the isomorphism
$\psi$ preserves both the anti-holomorphic multiplication and the level structures.  \qed

\section{Baily-Borel Compactification}\label{sec-BB}
Throughout this section we fix a square-free integer $d<0.$  Let $g \mapsto
\widehat{g} = N_{\beta} g N_{\beta}^{-1}$ be the resulting involution on
$\mathbf{Sp}(4,\mathbb R)$ with fixed point set $\phi(\mathbf{SL}(2,\mathbb C)).$  Let $\Gamma
\subset \mathbf{Sp}(4,\mathbb Z)$ be a torsion-free subgroup of finite index such that
$\widehat{\Gamma} = \Gamma$.  Let $\overline{\mathfrak h}_2$ be the partial Satake
compactification which is obtained by attaching rational boundary components of (complex)
dimension 0 and 1 (with the Satake topology).  The quotient $\overline{X} = \Gamma \backslash
\overline{\mathfrak h}_2$ is the Baily-Borel compactification of $X.$  It is a complex
projective algebraic variety.  Denote by
\[ \partial \overline{X} =\overline{X}-X = \partial_0\overline{X} \cup \partial_1\overline{X}
\]
the decomposition of the singular set into its (complex-) 0 and 1 dimensional strata.  The
involution $g \mapsto \hat g $ passes to an involution on $\overline{X}$ and defines a real
structure on $\overline{X}$, which we refer to as the ($\Gamma,\ {}^{\widehat{}}\ )$-real
structure.  Its fixed point set is the set $\overline{X}(\mathbb R)$ of real points of
$\overline{X}.$

Throughout the rest of \S \ref{sec-BB} we fix a level $N \ge 1,$ set $\Gamma = \Gamma_N$ of \S
\ref{subsec-arithmeticsubgroups}, and $X = \Gamma \backslash \mathfrak h_2.$

\begin{thm}\label{thm-noboundary}
Suppose $\mathcal O_d$ is principal.  If $d \equiv 1\ \mod 4,$ 
assume that $N$ is divisible by 4.  If $d=-2,$ assume that $N$ is even.  If $d=-1,$ assume
that $N$ is even and $N\ge 4.$  Then the 1-dimensional boundary strata of $\overline{X}$
contain no real points, that is, $\overline{X}(\mathbb R) \cap \partial_1 \overline{X} =
\phi.$ \end{thm}
In the next few sections we give separate proofs for $d \equiv 1\ \mod 4$, $d=-1$, and
$d=-2.$  Let $\mathbf{GSp}(4,\mathbb R)$ be the set of real matrices $g$ such that $gJ\tr{g} =
\lambda g$ for some $\lambda \in \mathbb R^{\times}$ (where $J = \left( \begin{smallmatrix} 0
& I \\ -I & 0 \end{smallmatrix} \right)$ is the standard symplectic form). It acts on
$\mathfrak h_2$ by fractional linear transformations and the center acts trivially.  Identify
the upper halfplane $\mathfrak h_1$ with the {\it standard} 1-dimensional boundary component
$F_1$ of $\mathfrak h_2$ by mapping $z = x+iy \in \mathfrak h_1$ to  $\left(
\begin{smallmatrix} z & 0\\ 0 & i\infty \end{smallmatrix} \right)$ (as a limit of $2\times 2$
complex matrices).  Let $P_1\subset \mathbf{GSp}(4,\mathbb R)$ denote the maximal parabolic
subgroup which normalizes this boundary component.  It acts on $\mathfrak h_1$ via the
projection  $\nu_h:P_1 \to \mathbf{GL}(2,\mathbb R)$ which is given by
\begin{equation}\label{eqn-projection}
\left(\begin{matrix} 
a&0&b&* \\ *&t&*&* \\ c&0&d&* \\ 0&0&0&s
\end{matrix}\right) \mapsto
\left(\begin{matrix} 
 a&0&b&0\\0&1&0&0\\c&0&d&0\\0&0&0&1
\end{matrix}\right) \leftrightarrow \left( \begin{matrix} a&b\\c&d \end{matrix} \right).
\end{equation}
For any $p \in P_1$ denote by $p_h = \nu_h(p).$  If $pJ\tr{p} = \lambda J$ then $\det (p_h) =
ad-bc=st= \lambda.$

Define $L_{\beta}: \mathfrak h_2 \to \mathfrak h_2$ by $L_{\beta}(y) = -\beta y
\tr{\beta}^{-1},$ cf. \S \ref{subsec-involutions}.  It is the fractional linear tansformation
corresponding to the matrix 
\begin{equation}\label{eqn-Lbeta}
 L_{\beta} = \left(\begin{smallmatrix} -\beta & 0 \\ 0 & \tr{\beta}
\end{smallmatrix} \right)\in \mathbf{GSp}(4,\mathbb Z).\end{equation}
(Note that $L_{\beta}J\tr{L}_{\beta} = -dJ.$)  
Denote by $\tilde Z = L_{\beta} \hat Z = L_{\beta}\cdot N_{\beta}\cdot \bar Z = - \bar Z$ for
$Z \in \mathfrak h_2.$  This is an anti-holomorphic involution which extends to the partial
compactification $\overline{\mathfrak h}_2$ and preserves the standard boundary component
$F_1.$  Let $g \mapsto \tilde g$ denotes the involution $\left( \begin{smallmatrix} A&B\\C&D
\end{smallmatrix}\right) \mapsto \left(\begin{smallmatrix} A & -B\\ -C & D \end{smallmatrix}
\right)$ on $\mathbf{Sp}(4,\mathbb R).$  Then $\widetilde{gZ} = \tilde g \tilde Z$ for all
$g\in \mathbf{Sp}(4,\mathbb R)$ and all $Z \in \mathfrak h_2.$  

Suppose $x'\in \partial_1\overline{X}(\mathbb R).$  Then there is a lift $x \in
\overline{\mathfrak h}_2$ which lies in some 1-dimensional rational boundary component, and
there exists $\gamma \in \Gamma$ such that $\hat x = \gamma x.$  Since $\mathbf{Sp}(4,\mathbb
Z)$ acts transitively on the set of all rational boundary components of a given dimension
there exists $g \in \mathbf{Sp}(4,\mathbb Z)$ and there exists $z \in F_1$ so that $x = gz.$ 
Then $\hat g \hat z = \hat x$ and $\tilde g \tilde z = \tilde x.$  It follows that
$pz = \tilde z=-\bar z\in F_1$ where
\begin{equation}\label{eqn-p} 
p= \tilde g^{-1} L_{\beta} \gamma g\in \mathbf{GSp}(4,\mathbb Z). \end{equation}
Hence $p \in P_1,$  $p$ is integral, and by (\ref{eqn-projection}), its Hermitian part 
\[ p_h= \left(\begin{matrix} p_{11}&p_{13} \\ p_{31}& p_{33} \end{matrix} \right) \in
\mathbf{GL}(2,\mathbb R)\] 
is also integral and has determinant $-d.$ 

\subsection{The case $d \equiv 1\ \mod 4$}\label{subsec-dequiv1}  Since $p\cdot z = -\bar z$
we find
\[ p_{11}z + p_{13} = -\bar z (p_{31}z + p_{33}) = -p_{31}z\bar z - p_{33} \bar z.\]
Comparing imaginary parts, $p_{11} = p_{33}$ so 
\begin{equation}\label{eqn-det}
-d = \det(p_h) = p_{11}^2 - p_{13}p_{31}.\end{equation}

If $g = \left( \begin{smallmatrix} A&B \\ C& D \end{smallmatrix} \right) \in
\mathbf{Sp}(4,\mathbb Z)$ then $\tilde g = g - \left( \begin{smallmatrix} 0 & 2B \\
2C & 0 \end{smallmatrix} \right).$  Therefore $g^{-1}\tilde g \equiv I \ \mod 2.$  Since
$\beta \equiv I\ \mod 2,$ the same holds for $L_{\beta}.$  Moreover $\gamma \in \Gamma_N
\subset \Gamma(N/2) \subset \Gamma(2)$ so
\[ p = (\tilde g^{-1} g) g^{-1}(L_{\beta} \gamma)g \equiv I\ \mod 2 \]
and the same holds for $p_h.$  In other words, $p_{11}$ is odd while $p_{13}, p_{31}$ are
even.  Therefore $p_{11}^2 - p_{13}p_{31} \equiv 1 \ \mod 4$ which contradicts
(\ref{eqn-det}) since $-d \equiv -1\ \mod 4.$  This completes the proof of Theorem
\ref{thm-noboundary} in the case $d \equiv 1\ \mod 4.$

\begin{lem}\label{lem-even} 
Suppose $d = -1$ or $d=-2.$   Let $g,\gamma \in \mathbf{Sp}_4(\mathbb Z)$ and suppose that
$\gamma\equiv I \ \mod 2.$  Suppose $p=\tilde g^{-1} L_{\beta} \gamma g \in P_1$, cf.
{\rm (\ref{eqn-p})}.  Denote the Hermitian part, cf. {\rm(\ref{eqn-projection})}, by
\begin{equation*}
p_h = \left( \begin{matrix} p_{11}&p_{13}\\p_{31}&p_{33} \end{matrix} \right).\end{equation*}
Then $p_{13}$ and $p_{31}$ are even.  If $d = -1$ then $p_h\in \text{\rm SL}_2(\mathbb Z)(2) $
lies in the principal congruence subgroup of level 2.
\end{lem}

\subsection{Proof}  Set $q = \tilde g^{-1} L_{\beta} g.$  Then $ p \equiv q\ \mod 2$ so $p_h
\equiv q_h\ \mod 2.$  If $ \tilde g^{-1} = \left(
\begin{smallmatrix}A&B\\C&D\end{smallmatrix} \right)$ then
\[
q = \left( \begin{matrix}
-A\beta\tr D + B\tr{\beta}\tr C & -A\beta\tr B + B\tr{\beta} \tr A \\
-C\beta\tr D + D\tr{\beta} \tr C & -C\beta\tr B + D\tr{\beta} \tr A
\end{matrix} \right) \]
Suppose $A = \left( \begin{smallmatrix} a_1&a_2\\a_3&a_4
\end{smallmatrix}\right)$ and similarly for $B$.  Direct computation gives
$q_{13}=-a_1b_2-da_2b_1+db_1a_2+b_2a_1=0$ and similarly $q_{31}=0.$ 
Hence $p_{13}$ and $p_{31}$ are even.  If $d = -1$ then $p_{11}p_{33}- p_{13}p_{31}=1$ which
implies that $p_{11}$ and $p_{33}$ are odd.  Hence $p_h \equiv I\ \mod 2.$  \qed

\subsection{Proof of Theorem \ref{thm-noboundary} in the case $d=-2$}
The same argument in \S \ref{subsec-dequiv1} leads to equation (\ref{eqn-det}):  $2=p_{11}^2 -
p_{13}p_{31}.$  By Lemma \ref{lem-even}, $p_{13}$ and $p_{31}$ are even.  This leads to a
contradiction whether $p_{11}$ is even or odd.  

\subsection{Proof of Theorem \ref{thm-noboundary} in the case $d=-1$}  
Return to (\ref{eqn-p}) and assume that $\gamma \in \Gamma_N \subset \Gamma(N)$ for some even
integer $N \ge 4.$  By Lemma \ref{lem-even}, $p_h \in \text{\rm SL}_2(\mathbb Z)(2).$  Since
$p_hz = -\bar z$ we can apply \cite{GT} \S 5 to find $Y\in \mathbb R_{+}$ and $h \in
\mathbf{SL}(2,\mathbb Z)$ such that $z = h\cdot iY.$  Using (\ref{eqn-projection}), the
element $h$ may also be regarded as lying in $P_1 \cap \mathbf{Sp}(4,\mathbb Z)$.  Define $g_1
= gh.$  Then the element
\[ v = (\tilde g_1)^{-1}L_{\beta} \gamma g_1 \in P_1 \cap \mathbf{Sp}(4,\mathbb Z)\]
fixes the point $iY \in \mathfrak h_1.$  Again by Lemma \ref{lem-even} its Hermitian part $v_h
\in \Gamma(2)$ (which is torsion-free).  Since it fixes $iY$ it must equal $\pm I.$  In
summary we may write
\[ v = \left( \begin{matrix} \alpha & * \\ 0 & \tr{\alpha}^{-1} \end{matrix} \right)\]
where 
\[ \alpha = \left( \begin{matrix} \epsilon_1 & 0 \\ * & \epsilon_2 \end{matrix} \right) \]
with $\epsilon_1 = \pm 1$ and $\epsilon_2 = \pm1.$  If $g_1 = \left( \begin{smallmatrix}
A&B\\C&D \end{smallmatrix} \right)$ with $A = \left( \begin{smallmatrix} a_{1} & a_{2} \\
a_3 & a_4 \end{smallmatrix} \right)$ (and similarly for $B$, $C$, and $D$), then $\tilde g_1 v
\equiv L_{\beta} g_1\ \mod N$ gives
\[ A \alpha \equiv -\beta A\ \mod N \text{ and } C\alpha \equiv \beta C\ \mod N.\]
From this we obtain $a_2 \equiv \epsilon_2 a_4 \equiv -a_2\ \mod N$ and $-c_2 \equiv
\epsilon_2 c_4 \equiv c_2 \ \mod N.$ Provided $N \ge 3$ this implies that the second
column of the matrix $g_1$ is $0\ \mod N$, which contradicts the fact that $\det(g_1) =1.$
This completes the proof of Theorem \ref{thm-noboundary} in the case $d=-1.$

\section{Rational Structure}\label{sec-rational}

As in \S \ref{subsec-numberfield}, fix a square-free integer $d<0$, let $\mathcal O_d$ be the
full ring of integers in $\mathbb Q(\sqrt{d})$ and let $g \mapsto \hat g = N_{\beta} g
N_{\beta}^{-1}$ be the resulting involution on $\mathbf{Sp}(4,\mathbb R)$ and 
$Z \mapsto \hat Z = \beta \overline{Z} \tr{\beta}^{-1}$ the resulting  
anti-holomorphic involution  on $\mathfrak h_2.$  In this section we make no further
assumptions on $d.$ Fix a level $N \ge 1$ and, as in \S \ref{subsec-arithmeticsubgroups} let
$ \Gamma = \Gamma_N =  \Gamma(N) \cap \widehat{\Gamma}(N).$
\begin{thm}  The Baily-Borel compactification $\overline{X}$ admits the
structure of a complex projective variety which is defined over $\mathbb Q$ such that the
resulting real structure, when restricted to $X$ agrees with the
real structure of Theorem \ref{thm-realpoints}. \end{thm}

The proof will occupy the rest of this section.

\begin{prop} The complex vectorspace of (holomorphic) $\Gamma$-modular forms on $\mathfrak
h_2$ is 
spanned by modular forms with rational Fourier coefficients. \end{prop}
\subsection{Proof}  This follows directly from \cite{Shimura} who proves the following in
Theorem 3 (ii).  Let $\mathbf G = \mathbf{GSp}_4$, let $\mathbb A$ be the adeles of $\mathbb
Q$ and let $S \subset \mathbf G(\mathbb A)^{+}$ be an open subgroup containing $\mathbb
Q^{\times}\mathbf G(\mathbb R)^{+}$ such that $S/ \mathbb Q^{\times} \mathbf G(\mathbb R)^{+}$
is compact.  (Here, $+$ denotes the identity component.)  Let $\Gamma = S \cap \mathbf
G(\mathbb Q).$  Suppose that
\[ \Delta = \left\{ \left( \begin{matrix} I & 0 \\ 0 & tI \end{matrix} \right)
\biggl|\biggr. \ t \in \underset{p}{\Pi}\ \mathbb Z_p^{\times} \right\} \subset S.\]
Then the complex vectorspace of $\Gamma$-modular forms with weight $k$ on $\mathfrak h_2$ is
spanned by those forms whose Fourier coefficients are in the field $k_S = \mathbb Q.$  To
apply this to our setting, define
\begin{align*}
S(N) &= \left\{ x \in \mathbf G(\mathbb A)^{+} \left|\right.  \ x_p \equiv \left(
\begin{smallmatrix}
I & 0 \\ 0 & a_pI \end{smallmatrix} \right) \mod{N \cdot \mathbb Z_p}, \exists a_p \in
\mathbb Z_p^{\times} \right\} \cdot \mathbb Q^{\times} \\
S^{\beta}(N) &= N_{\beta} S(N) N_{\beta}^{-1}.
\end{align*}
(Here, $x_p$ denotes the $p$-component of $x.$)
It is easy to see that each of these contains $\Delta$, hence Shimura's hypothesis is
satisfied.  If $S = S(N) \cap S^{\beta}(N)$ then $\Gamma_S = S \cap \mathbf G(\mathbb Q) =
(\Gamma(N) \cap \widehat{\Gamma}(N)) \cdot \mathbb Q^{\times}.$ \qed

\subsection{}  Let $I_- = \left( \begin{smallmatrix} I & 0 \\ 0 & -I \end{smallmatrix}
\right).$  Its action by fractional linear transformations maps the Siegel lower halfspace
$\mathfrak h^{-}_2$ to the upper halfspace $\mathfrak h_2$, that is, $I_- \cdot Z = -Z.$
Hence,
for any holomorphic mapping $f:\mathfrak h_2 \to \mathbb C$ we may define $f^{\beta}:\mathfrak
h_2 \to \mathbb C$ by 
\begin{equation*}  f^{\beta}(Z) = f(I_- \cdot N_{\beta} \cdot Z) = f(-\beta Z
\tr{\beta}^{-1}). \end{equation*}

\begin{prop}  If $f:\mathfrak h_2 \to \mathbb C$ is a holomorphic $\Gamma$-modular form of
weight $k$ then $f^{\beta}$ is also a holomorphic $\Gamma$-modular form of weight $k$, and
\begin{equation}\label{eqn-modular}
f^{\beta}(\widehat{Z}) = \overline{f^{\beta}(Z)}\end{equation}
for all $Z \in \mathfrak h_2.$ \end{prop}

\subsection{Proof}
Suppose that $f(\gamma \cdot Z) = j(\gamma,Z)^kf(Z)$ for all $\gamma \in \Gamma$ and all $Z
\in\mathfrak h_2$ where 
\begin{equation*}
j\left( \left( \begin{smallmatrix} A & B \\ C & D \end{smallmatrix} \right), Z \right)
= \det(CZ+D) \end{equation*}
is the standard automorphy factor.  Then $j(I_-N_{\beta},Z) = \det(-\tr{\beta})$ is
independent of $Z$.  Let $\gamma \in \Gamma$ and set
\[ \gamma ' = I_- N_{\beta} \gamma N_{\beta}^{-1} I_-^{-1} \in \Gamma.\]  Then
\begin{align*}
f^{\beta}(\gamma \cdot Z) &= f(\gamma' \cdot I_- N_{\beta} \cdot Z)\\
&= j(\gamma', I_-N_{\beta}\cdot Z)^k f^{\beta}(Z)\\
&= \det(-\tr{\beta})^k j(\gamma,Z)^k \det(-\tr{\beta})^{-k} f^{\beta}(Z)\\
&= j(\gamma,Z)^kf^{\beta}(Z).
\end{align*}
which shows that $f^{\beta}$ is $\Gamma$-modular of weight $k$.  Next, with respect to the
standard maximal parabolic subgroup $P_0$ (which normalizes the standard 0-dimensional
boundary component), the modular form $f$ has a Fourier expansion,
\begin{equation*}
f(Z) = \sum_s a_s \exp\left({2\pi i \langle s,Z \rangle}\right)
\end{equation*}
which is a sum over lattice points $s\in L^*$ where $L = \Gamma \cap Z(\mathcal U_0)$ is the
intersection of $\Gamma$ with the center of the unipotent radical $\mathcal U_0$ of $P_0$ and
where $a_s \in \mathbb Q.$  Then
\begin{align*}
f(\hat Z)   &= \sum_s a_s \exp\left( 2 \pi i \langle s, \beta \overline{Z }\tr{\beta}^{-1}
\rangle\right) \\
&= \overline{\sum_s a_s  \exp \left( 2 \pi i \langle s, -\beta Z \tr{\beta}^{-1} \rangle
\right)} \\
&= \overline{f(-\beta Z \tr{\beta}^{-1})} = \overline{f^{\beta}(Z)}.  \qed \end{align*}

\subsection{}

The Baily-Borel compactification $\overline{X}$ of $X$ is the obtained by
embedding $X$ holomorphically into $\mathbb C \mathbb P^m$ using $m+1$ ($\Gamma$-)modular
forms (say $f_0, f_1,\ldots, f_m$) of some sufficiently high weight $k$,  and then
taking the closure of the image.  Define an embedding $\Phi: X \to \mathbb C
\mathbb P^{2m+1}$ by
\begin{equation*}
\Phi(Z) = \left(f_0(Z):f_1(Z):\ldots :f_m(Z): f^{\beta}_0(Z):f^{\beta}_1(Z)\ldots :
f^{\beta}_m(Z)\right).
\end{equation*}
Denote these homogeneous coordinate functions by $x_j = f_j(Z)$ and $y_j =
f^{\beta}_j(Z).$  Define an involution $\sigma: \mathbb C \mathbb P^{2m+1} \to \mathbb
C \mathbb P^{2m+1}$ by $\sigma(x_j) = \bar y_j$ and $\sigma(y_j) = \bar x_j.$
Then equation (\ref{eqn-modular}) says that this involution is compatible with the
embedding $\Phi$, that is, for all $Z \in X$ we have:
\begin{equation*}
\sigma\Phi(Z) = \Phi(\widehat Z).
\end{equation*}  

Define $\Psi: \mathbb C \mathbb P^{2m+1} \to \mathbb C \mathbb P^{2m+1}$ by
setting $\xi_j = x_j+y_j$ and $\eta_j = i(x_j - y_j)$ for $0 \le j \le m.$  
Let $Y = \Psi \Phi(X)$ and let $\overline{Y}$ denote its closure.
\begin{prop}
The composition $\Psi \Phi:X \to \mathbb C \mathbb P^{2m+1}$ is a holomorphic
embedding which induces an isomorphism of complex algebraic varieties $\overline{X} \to
\overline{Y}.$  The variety $\overline{Y}$ is defined over the rational numbers, and the real
points of $Y$ are precisely the image of those points $Z \in X$ such that $\hat Z = Z.$
\end{prop} 
\subsection{Proof.}  The image $\Psi\Phi(X)$ is an algebraic subvariety of
projective space, which is preserved by complex conjugation so it is defined
over $\mathbb R.$  The real points are obtained by setting $\bar \xi_j
=\xi_j$ and $\bar \eta_j = \eta_j$ which gives $\bar x_j = y_j$ and $\bar y_j =
x_j$ hence $\Phi(Z) = \sigma\Phi(Z),$ or $Z = \hat Z.$ The Fourier
coefficients of $\xi_j$ and $\eta_j$ are in $\mathbb Q[i]$ so the image
$\Psi\Phi(X)$ is defined over $\mathbb Q[i].$  Since it is also invariant under
$\text{Gal}(\mathbb C/\mathbb R)$, it follows that 
$\Psi\Phi(X)$ is defined over $\mathbb Q.$  \qed

\subsection{Remark}  The usual embedding $(f_0:f_1:\ldots:f_m):X \hookrightarrow \mathbb C
\mathbb P^m$ determines the usual rational structure on $\overline{X}$, and the resulting
complex conjugation is that induced by $Z \mapsto -\overline{Z}$ for $Z \in\mathfrak h_2.$


\begin{thebibliography}{AMRT}


\bibitem[C]{Comessatti} H. Comessatti, Sulle variet\`a abeliane reali I, II.  Ann. Mat. Pura.
Appl. {\bf 2} (1924) pp. 67-106 and {\bf 4} (1926) pp. 27-72.


\bibitem[GT]{GT} M. Goresky and Y.-S. Tai, The moduli space of real abelian varieties with
level structure.

\bibitem[La]{Lange} H. Lange and C. Birkenhake, {\bf Complex Abelian
Varieties}, Grundlehren Math. Wiss. {\bf 302}, Springer Verlag, Berlin,
1992.

\bibitem[N]{Nygaard} N. Nygaard, Cohomology of Siegel 3-folds.  Comp. Math. {\bf 97} (1995),
173-186.

\bibitem[Si]{Silhol} R.~Silhol, Compactifications of moduli spaces in real algebraic geometry.
Inv. Math. {\bf 107} (1992), pp. 151-202.


\bibitem[Si2]{Silhol2} R.~Silhol, Real abelian varieties and the theory of Comessatti.  Math.
Z. {\bf 181} (1982), pp. 345-362.

\bibitem[Sh]{Shimura} G.~Shimura, On the Fourier coefficients of modular forms of several
variables.  Nachr. Akad. Wiss. G\"ottingen Math-Phys.  {\bf 17} (1975), 261-268.

\bibitem[St]{Stark} H.~Stark, A complete determination of the complex quadratic fields of
class-number one.  Mich. Math. J. {\bf 14} (1967), 1-27.


\end{thebibliography}
\end{document}